\documentclass{article}
\usepackage{fullpage,amssymb,amsmath,amsthm,amsfonts}

\def\polylog{\mathop{\rm polylog}\nolimits}

\newcommand{\bmat}{\left[ \begin{array}}
\newcommand{\emat}{\end{array} \right]}
\newcommand{\diag}{{\rm diag}}
\newcommand{\sign}{{\rm sign}}
\newcommand{\sep}{{\rm sep}}
\newcommand{\VEC}{{\rm vec}}

\newcommand{\mtf}{\lfloor \frac{m}{2} \rfloor}
\newcommand{\mtc}{\lceil \frac{m}{2} \rceil}
\newcommand{\mt}{\frac{m}{2}}

\newcommand{\ignore}[1]{}

\newtheorem{theorem}{Theorem}[section]

\newtheorem{lemma}[theorem]{Lemma}
\newtheorem{corollary}[theorem]{Corollary}

\theoremstyle{definition}
\newtheorem{definition}[theorem]{Definition}

\title{Fast Linear Algebra is Stable}
\author{James Demmel\thanks{Mathematics Department and CS Division,
University of California, Berkeley, CA 94720. 
The author acknowledges support of NSF under grants CCF-0444486, ACI-00090127, 
CNS-0325873 and of DOE under grant DE-FC02-01ER25478.},
Ioana Dumitriu\thanks{Mathematics Department, University of Washington, Seattle, WA 98195. },
and
Olga Holtz\thanks{Mathematics Department, University of California, Berkeley,
CA 94720.},
}

\date{August 28, 2007}
\begin{document}
\maketitle

\begin{abstract}
In \cite{FastMatMulIsStable} we showed that a large class 
of fast recursive matrix multiplication algorithms is stable 
in a normwise sense, and that in fact if multiplication of 
$n$-by-$n$ matrices can be done by any algorithm
in $O(n^{\omega + \eta})$ operations for any $\eta > 0$, 
then it can be done stably in $O(n^{\omega + \eta})$
operations for any $\eta > 0$. Here we extend this result to show that essentially
all standard linear algebra operations, including LU decomposition, QR decomposition, 
linear equation solving, matrix inversion, solving least squares problems,
(generalized) eigenvalue problems and the singular value decomposition
can also be done stably (in a normwise sense) in $O(n^{\omega + \eta})$ operations.
\end{abstract}

\section{Introduction}
\label{sec_Intro}

Matrix multiplication is one of the most fundamental operations in numerical linear 
algebra. Its importance is magnified by the number of other problems (e.g., computing 
determinants, solving systems of equations,  matrix inversion, LU decomposition, 
QR decomposition, least squares problems etc.) 
that are reducible to it \cite{burgisserclausenshokrollahi,Heller78,BorodinMunro75}. 
This means that an algorithm for multiplying $n$-by-$n$ matrices in $O(n^{\omega})$ 
operations can be converted into an algorithm for these other linear algebra operations
that also runs in $O(n^{\omega})$ operations.

In this paper we extend this result to show that if the matrix multiplication 
algorithm is stable in a normwise sense discussed below, then essentially {\em all}
linear algebra operations can also be done stably, in time $O(n^{\omega})$
or $O(n^{\omega + \eta})$, for arbitrary $\eta > 0$.
For simplicity, whenever an exponent contains ``$+ \eta$'', we
will henceforth mean ``for any $\eta > 0$.''

In prior results~\cite{FastMatMulIsStable} we showed that any fast matrix 
multiplication algorithm running in time $O(n^{\omega + \eta})$ 
was either stable or could be converted into a stable algorithm that also 
ran in $O(n^{\omega + \eta})$ operations. 
Combined with the results in this paper, this lets us state
that all linear algebra operations can also be done stably in 
$O(n^{\omega + \eta})$ operations. 

More precisely,
some of our results (see Theorem~\ref{thm_1} in Section~\ref{sec_LogStable})
may be roughly summarized by saying that $n$-by-$n$ matrices
can be multiplied in $O(n^{\omega + \eta})$ operations 
{\em if and only if} $n$-by-$n$ matrices can be inverted stably 
in $O(n^{\omega + \eta})$ operations.
We need to use a little bit of extra precision to make this claim,
and count operations carefully;
the cost of extra precision is accounted for by the $O(n^{\eta})$ factor.

Other results (see Section~\ref{sec_Simultaneous})
may be summarized by saying that if $n$-by-$n$ matrices
can be multiplied in $O(n^{\omega + \eta})$ {\em arithmetic} operations, 
then we can compute the QR decomposition stably 
(and so solve linear systems and least squares problems stably)
in $O(n^{\omega + \eta})$ {\em arithmetic} operations. 
These results do not require extra precision, which is why we only 
need to count arithmetic operations, not bit operations.

The QR decomposition will then be used to stably compute a rank-revealing
decomposition, 
compute the (generalized) Schur form, 
and compute the singular value decomposition, 
all in $O(n^{\omega + \eta})$ {\em arithmetic} operations. 
To compute (generalized) eigenvectors from the Schur form 
we rely on solving the (generalized) Sylvester equation 
all of which can be done stably in $O(n^{\omega + \eta})$ {\em bit} operations.



Now we become more precise about our notions of stability. We say an algorithm
for multiplying $n$-by-$n$ square matrices $C = A \cdot B$ is {\em stable}
if the computed result $C_{comp}$ satisfies the following normwise error bound:
\begin{equation} 
\| C_{comp}-C\|\leq \mu(n) \varepsilon \|A\| \, \|B \| +O(\varepsilon^2), \label{gen_bound} 
\end{equation}
where $\varepsilon$ is machine epsilon (bounds the roundoff error) and
$\mu(n)$ is a (low degree) polynomial, i.e., $\mu(n)=O(n^c)$ for some constant $c$. 
Note that one can easily switch from one norm to another at the expense of picking up 
additional factors that will depend on $n$, using the equivalence of norms on a 
finite-dimensional space, thereby changing the constant $c$ slightly.
The bound~(\ref{gen_bound}) was first obtained for Strassen's $O(n^{2.81})$ 
algorithm~\cite{Strassen} by Brent (\cite{brent3,higham90b}, \cite[chap. 23]{higham96})
and extended by Bini and Lotti~\cite{BiniLotti} to a larger
class of algorithms. In prior work~\cite{FastMatMulIsStable} we
showed that such a bound holds for a new class of fast algorithms
depending on group-theoretic methods~\cite{CU} and~\cite{CKSU},
which include an algorithm that runs asymptotically as fast as the
fastest known method due to Coppersmith and Winograd~\cite{CoppersmithWinograd},
which runs in about $O(n^{2.38})$ operations.
Using a result of Raz~\cite{Raz}, that work also showed that any
fast matrix multiplication algorithm running in $O(n^{\omega + \eta})$ arithmetic
operations 
can be converted to one that satisfies (\ref{gen_bound})
and also runs in $O(n^{\omega + \eta})$ arithmetic operations.

In Section~\ref{sec_Block} we begin by reviewing conventional 
block algorithms used in practice in libraries like
LAPACK \cite{lapackmanual3} and ScaLAPACK \cite{scalapackmanual}.
The normwise backward stability of these algorithms was demonstrated in 
\cite{DemmelHigham92,higham90b,demmelhighamschreiber,higham96}
using (\ref{gen_bound}) as an assumption; this means that these
algorithms are guaranteed to produce the exact answer (e.g., solution
of a linear system) for a matrix $\hat{C}$ close to the actual 
input matrix $C$, where close means close in norm: 
$\|\hat{C} - C \| = O(\varepsilon) \|C\|$. Here the $O(\varepsilon)$ is
interpreted to include a factor $n^c$ for a modest constant $c$.

What was not analyzed in this earlier work was the speed of these
block algorithms, assuming fast matrix multiplication. 
In Section~\ref{sec_Block} we show that the optimal choice of
block size lets these block algorithms run only as fast as 
$O(n^{\frac{9 - 2\gamma}{4 - \gamma}})$ operations, 
where $O(n^{\gamma})$ is the operation count of matrix multiplication.
(We use $\gamma$ instead of $\omega + \eta$ to simplify notation.)
Even if $\gamma$ were to drop from 3 to 2, the exponent $\frac{9 - 2\gamma}{4 - \gamma}$
would only drop from 3 to 2.5. While this is an improvement, we shall do better.

In Section~\ref{sec_LogStable} we consider known divide-and-conquer 
algorithms for reducing the complexity of matrix inversion
to the complexity of matrix multiplication. These algorithms are 
not backward stable in the conventional sense. However, 
we show that they can achieve the same
forward error bound (bound on the norm of the error in the output) as
a conventional backward stable algorithm, provided that they use
just $O(\log^p n)$ times as many bits of precision in each arithmetic
operation (for some $p>0$) as a conventional algorithm. We call such algorithms
{\em logarithmically stable.\/}  Incorporating the cost of this
extra precise arithmetic in the analysis only increases the total cost
by a factor at most $\log^{2p} n$. Thus, if there are matrix multiplication
algorithms running in $O(n^{\omega + \eta})$ operations for any $\eta > 0$,
then these logarithmically stable algorithms for operations like
matrix inversion also run in $O(n^{\omega + \eta})$ operations 
for any $\eta > 0$, and achieve the same error bound as a conventional algorithm.

In Section~\ref{sec_QRStable} we analyze a divide-and-conquer algorithm 
for QR decomposition described in~\cite{ElmrothGustavson2000} that is simultaneously 
backward stable in the conventional normwise sense (i.e. without extra precision), 
and runs in $O(n^{\omega + \eta})$
operations for any $\eta > 0$. This may be in turn used to solve linear
systems, least squares problems, and compute determinants equally stably and fast.
We apply the same idea to LU decomposition in Section~\ref{sec_LUStable}
but stability depends on a
pivoting assumption similar to, but slightly stronger than, the usual
assumption about the stability of partial pivoting.

In Section~\ref{sec_URV} we use the QR decomposition to compute a rank revealing
$URV$ decomposition of a matrix $A$. This means that $U$ and $V$ are orthogonal,
$R$ is upper triangular, and $R$ reveals the rank of $A$ in the following sense:
Suppose
$\sigma_1 \geq \cdots \geq \sigma_n$ are the singular values of $A$. Then for each $r$,
$\sigma_{\min} (R(1:r,1:r))$ is an approximation of $\sigma_r$
and
$\sigma_{\max} (R(r+1:n,r+1:n))$ is an approximation of $\sigma_{r+1}$.
(Note that if $R$ were diagonal, then the URV decomposition would be
identical to the singular value decomposition, and these approximations 
would be exact.) Our algorithm will be {\em randomized}, in the sense that 
the approximations of $\sigma_r$ and $\sigma_{r+1}$ are reasonably accurate 
with high probability.

In Section~\ref{sec_Schur}, we use the QR and URV decompositions
in algorithms for the (generalized) Schur form of nonsymmetric matrices
(or pencils) \cite{baidemmelgu94}, lowering their complexity to $O(n^{\omega + \eta})$
arithmetic operations while maintaining normwise backward stability.
The singular value decomposition may in turn be reduced to solving 
an eigenvalue problem with the same complexity (Section~\ref{sec_SymEigSVD}).
Computing (generalized) eigenvectors can only be done in a logarithmically stable
way from the (generalized) Schur form.
We do this by providing a logarithmically stable algorithm for solving the
(generalized) Sylvester equation, and using this to compute eigenvectors.
A limitation of our approach is that to compute all the eigenvectors
in $O(n^{\omega + \eta})$ bit operations, all the eigenvectors may 
in the worst case have a common error bound that depends on the
worst conditioned eigenvector.

\section{Conventional Block Algorithms}
\label{sec_Block}

A variety of ``block algorithms'' that perform most of their operations
in matrix multiplication are used in practice 
\cite{lapackmanual3,scalapackmanual}
and have been analyzed in the literature
\cite{DemmelHigham92,higham90b,demmelhighamschreiber,higham96}, 
and it is natural to consider these conventional algorithms first. 
For example, \cite{DemmelHigham92} does a general error
analysis of block algorithms for LU factorization, QR factorization,
and a variety of eigenvalue algorithms using the bound (\ref{gen_bound}),
and shows they are about as stable as their conventional counterparts.
What was not analyzed in \cite{DemmelHigham92} was the complexity,
assuming fast matrix multiplication.

We will use the notation $MM(p,q,r)$ to mean
the number of operations to multiply a $p$-by-$q$ times 
a $q$-by-$r$ matrix; when $q \leq p, r$, this is done by
$\frac{p}{q} \cdot \frac{r}{q}$ multiplications of $q$-by-$q$ matrices,
each of which costs $MM(q,q,q) = O(q^{\gamma})$ for some $2 < \gamma \leq 3$.
(We use $\gamma$ instead of $\omega + \eta$ in order to simplify notation.)
Thus $MM(p,q,r) = O(pq^{\gamma - 2}r)$.
Similarly, when $p$ is smallest $MM(p,q,r) = O(p^{\gamma - 2}qr)$, and so on.
Also we will abbreviate $MM(n,n,n) = MM(n)$.

Consider block LU factorization with pivoting. 
Given a choice of block size $b$,
the algorithm breaks the $n$-by-$n$ matrix $A$ into blocks of $b$ columns,
then LU factorizes each such block using the conventional algorithm,
and then updates the trailing part of the matrix using fast matrix
multiplication. This may be expressed as
\begin{equation}
A = \bmat{cc} A_{11} & A_{12} \\ A_{21} & A_{22} \emat
 = P \cdot \bmat{cc} L_{11} & 0 \\ L_{21} & I \emat \cdot
           \bmat{cc} U_{11} & U_{12} \\ 0 & \hat{A}_{22} \emat
\end{equation}
where $A_{11}$ and $L_{11}$ are $b$-by-$b$, and $P$ is a permutation.
Thus the steps of the algorithm are as follows:
\begin{description}
\item[a.] Factor $\bmat{c} A_{11} \\ A_{21} \emat = P \bmat{c} L_{11} \\ L_{21} \emat U_{11}$. 
\item[b.] Apply $P^T$ to the rest of the matrix columns (no arithmetic operations).
\item[c.] Solve the triangular system $L_{11} U_{12} = A_{12}$ for $U_{12}$.
\item[d.] Update the Schur complement $\hat{A}_{22} = A_{22} - L_{21} U_{12}$.
\item[e.] Repeat the procedure on $\hat{A}_{22}$.
\end{description}
The conventional algorithm  \cite{demmelMA221,golubvanloan}
for steps (a) and (c) costs $O(nb^2)$.
Step (d) involves matrix multiplication at a cost 
$MM(n-b,b,n-b) = O(n^2b^{\gamma-2})$. 
Repeating these steps $n/b$ times makes the total cost
$O(n^2 b + n^3 b^{\gamma - 3})$.


To roughly minimize this cost, choose $b$ to make $n^2b = n^3 b^{\gamma - 3}$,
yielding $b = n^{\frac{1}{4 - \gamma}}$ and
$\# ops = O(n^{\frac{9-2\gamma}{4-\gamma}})$.
For $\gamma \approx 3$, the cost is near the usual $O(n^3)$, but
as $\gamma$ decreases toward 2, $b$ drops to $n^{1/2}$ but the $\# ops$
only drops to $O(n^{2.5})$. 

This same big-O analysis applies to $QR$ factorization:
When $A$ is $n$-by-$m$ and real and $n \geq m$, then we can 
write $A=QR$ where $Q$ is $n$-by-$n$ and orthogonal and $R$ 
is $n$-by-$m$ and upper triangular. 
We will represent $Q$ compactly 
by using the $WY$ representation of $Q$
\cite{BischofVanLoan87}: 
$Q^T$ can be written $Q^T = I - WY$,
where $W$ and $Y^T$ are both $n$-by-$m$ and lower triangular,
$W$'s columns all have 2-norm equal to 1, and
$Y$'s columns all have 2-norm equal to 2.
(An even lower-memory version of this algorithm 
\cite{Schreiber87a}
is used in practice \cite{lapackmanual3,scalapackmanual},
but we use \cite{BischofVanLoan87} for simplicity
of presentation.)
The conventional algorithm using the $WY$ representation
or variations costs $O(nm^2)$ operations to compute $R$, 
$O(nm^2)$ operations to compute $W$ and $Y$, 
and $O(n^2m)$ operations to explicitly construct $Q$, or
$O(n^3)$ in the square case \cite{BischofVanLoan87}.

The algorithm
for block QR factorization is entirely analogous to block
LU factorization, processing the matrix in blocks of $b$
columns at a time, updating the trailing part of the matrix
using fast matrix multiplication, based on the following
identity, where $A_1$ is $n$-by-$b$ and $A_2$ is $n$-by-$(m-b)$:
\begin{eqnarray}
\label{eqn_QR}
A & = & [ A_1 , A_2 ]  =  [ Q_1 R_1 , A_2 ] 
                   =  Q_1 [ R_1 , Q_1^T A_2 ]  \nonumber \\
                   & = & Q_1 [ R_1, (I - W_1 Y_1)A_2]
                   =  Q_1 [ R_1 , A_2 - W_1 ( Y_1 A_2) ] 
                   =  Q_1 [ R_1 , \hat{A}_2 ]
\end{eqnarray}
where $Q_1^T = I - W_1 Y_1$.
The cost of this step is $O(nb^2)$ for $A_1 = Q_1 R_1$ plus
\[
MM(b,n,m-b) + MM(n,b,m-b) + n(m-b) = O(nm b^{\gamma-2})
\]
for $\hat{A}_2$.
Repeating this procedure $(m-b)/b$ times on the last $n-b$ rows of 
$\hat{A}_2 = \bmat{c} \hat{A}_{21} \\ \hat{A}_{22} \emat$ eventually yields
$\hat{A}_{22} = Q_2 R_2 = (I - W_2 Y_2)^T R_2 $. 
Combining this with (\ref{eqn_QR}) yields
\begin{equation}
A = Q_1 \cdot \bmat{cc} R_{11} & \hat{A}_{12} \\ 0 & Q_2 R_2 \emat
    = Q_1 \cdot \bmat{cc} I & 0 \\ 0 & Q_2 \emat \cdot 
      \bmat{cc} R_{11} & \hat{A}_{12} \\ 0 & R_2 \emat 
    \equiv Q_1 \cdot \hat{Q}_2 \cdot R
    \equiv  Q \cdot R
\end{equation}
In practice we leave $Q = Q_1 \cdot \hat{Q}_2$ in this factored form
(note the $Q_2$ will also be a product of factors)
since that is faster for subsequent purposes, like solving least
squares problems. Thus the cost is bounded by $(m/b)$ times the cost of
(\ref{eqn_QR}), namely $O(nmb + nm^2 b^{\gamma-3})$. 
When $n=m$, this is the same cost as for block Gaussian elimination.
In the general case of $m \leq n$, we again roughly minimize
the cost by choosing $b$ so $nmb = nm^2 b^{\gamma-3}$, namely
$b = m^{1/(4-\gamma)}$, leading to a cost of 
$O(nm^{\frac{5-\gamma}{4-\gamma}})$. As $\gamma$ drops from
3 toward 2, this cost drops from $O(nm^2)$ toward $O(nm^{1.5})$.

If we wish, we may also multiply out the $Q_i^T$ factors into a 
matrix of the form $I - WY$ where $W$ and $Y^T$ are $n$-by-$m$
and lower triangular. Getting $W$ and $Y$ costs $O(nm^2b^{\gamma-3})$, 
and multiplying out $I - WY$ costs an additional
$O(n^2 m b^{\gamma - 3})$. This does not change
the cost in a big-O sense when $m/n$ is bounded below.
The following equation shows how:
\begin{eqnarray*}
Q^T & = & (I - \bmat{c} 0 \\ W_2 \emat \cdot [0, Y_2]) \cdot (I-W_1 \cdot Y_1) \\
      & \equiv & (I - \hat{W}_2 \cdot \hat{Y}_2) \cdot (I-W_1 \cdot Y_1) \\
      & = & (I - [\hat{Q}_2^T \cdot W_1, \hat{W}_2] \cdot [Y_1; \hat{Y}_2]) \\
      & = & (I - [W_1 - \hat{W}_2 \cdot (\hat{Y}_2 \cdot W_1) , \hat{W}_2] \cdot [Y_1; \hat{Y}_2]) \\
      & \equiv & I - WY
\end{eqnarray*}
(here we have used Matlab notation like $[Y_1; \hat{Y}_2]$ to stack
$Y_1$ on top of $\hat{Y}_2$).
Now the cost minimizing $b$ leads to a cost of 
$O(n^{\frac{5-\gamma}{4-\gamma}}m)$.

In summary, conventional block algorithm guarantee stability
but can only reduce the operation count to $O(n^{2.5})$ even
when matrix multiplication costs only $O(n^2)$. 
To go faster,  other algorithms are needed.

\ignore{

\begin{tabbing}
jnk \= jnk \= jnk \= jnk \= jnk \kill
Block $QR$ Factorization \\
\> function $[W, Y, R] = BlockQR(A,n)$ \\
\> ... Factor $A = (I+WY)^TR$ where $Q^T = I+WY$ is orthogonal and $R$ is upper triangular, \\
\> \> for $j=1$ to $n/b$ ... assume $b$ divides $n$ \\
\> \> \> $i = (j-1)b + 1$, $end = i+b-1$ \\
(1) \> \> \> Factor $A(i:n, i:end) = (I - \hat{W} \cdot \hat{Y})^T \cdot \hat{R}$ \\
\> \> \> ... $b$ columns of $A$ using conventional QR factorization, where 
             $\hat{W}$ and $\hat{Y}^T$ are $n-i+1$-by-$b$ lower trapezoidal, \\
\> \> \> ... and $\hat{R}$ is $n-1+1$-by-$b$ upper triangular \\
\> \> \> $L(i:n, i:end) = \hat{L}$, $U(i:end, i:end) = \hat{U}$, 
            update $P$ with $\hat{P}$ \\
(2) \> \> \> Compute $U(i:end, end+1:n) = \hat{L}^{-1} A(i:end, end+1:n) \\
\> \> \> ... $b$ rows of $U$ column by column by conventional substitution \\
(3) \> \> \> Update $A(end+1:n, end+1:n) = A(end+1:n, end+1:n) - L(end+1:n, i:end) 
             \cdot $U(i:end,end+1:n)$ \\

The boundedness of $W$ and $Y$ can
linear algebra algorithms as well. To go faster than this requires
different algorithms.

} 

\section{Logarithmically Stable Algorithms}
\label{sec_LogStable}

Our next class of fast and stable algorithms will abandon
the strict backward stability obtained by conventional algorithms
or their blocked counterparts in the last section in order
to go as fast as matrix multiplication. Instead, they will use
extra precision in order to attain roughly the same forward
errors as their backward stable counterparts. We will show
that the amount of extra precision is quite modest, and grows
only proportionally to $\log n$. Depending on exactly how
arithmetic is implemented, this will increase the cost of the
algorithm by only a $\polylog(n)$ factor, i.e. a polynomial in
$\log n$. For example,
if matrix multiplication costs $O(n^{\gamma})$ with $2 < \gamma \leq 3$,
then for a cost of $O(n^{\gamma} polylog(n)) = O(n^{\gamma + \eta})$
for arbitrarily tiny $\eta > 0$ one can invert matrices as
accurately as a backward stable algorithm. We therefore call
these algorithms {\em logarithmically stable}.

To define logarithmic stability more carefully,
suppose we are computing $y = f(x)$. Here $x$ could denote
a scalar, matrix, or set of such objects, equipped with
an appropriate norm. For example,
$y = f(\{A,b\}) = A^{-1}b$ is the solution of $Ay=b$.
Let $\kappa_f(x)$ denote the condition number of $f()$,
i.e. the smallest scalar such that
\[
\frac{\|f(x + \delta x) - f(x) \|}{\|f(x)\|} \leq
\kappa_f(x) \cdot \frac{\| \delta x \|}{\|x\|} + 
O( \left( \frac{\| \delta x\|}{\|x\|} \right)^2 ) \; \; .
\]
Let $alg(x)$ be the result of a backward stable algorithm for $f(x)$,
i.e. $alg(x) = f(x + \delta x)$ where $\|\delta x\| = O( \varepsilon ) \| x \|$.
This means the relative error in $alg(x)$ is bounded by
$$\frac{\|alg(x) - f(x)\|}{\|f(x)\|} = O( \varepsilon ) \kappa_f(x) + O( \varepsilon^2 ).$$

\begin{definition}[Logarithmic Stability]
Let $alg_{ls}(x)$ be an algorithm for $f(x)$, where the
``size'' (e.g., dimension) of $x$ is $n$.
If the relative error in $alg_{ls}(x)$ is bounded by 
\begin{equation}
\label{eqn_LogStable}
\frac{\|alg_{ls}(x) - f(x)\|}{\|f(x)\|} = O( \varepsilon ) \kappa_f^{\chi(n)}(x) 
+ O( \varepsilon^2 )
\end{equation}
where $\chi(n) \geq 1$ is bounded by a polynomial in $\log n$, then we
say $alg_{ls}(x)$ is a {\em logarithmically stable} algorithm for $f(x)$.
\end{definition}

\begin{lemma}
\label{lemma_LogStable}
Suppose $alg_{ls}(x)$ is a logarithmically stable algorithm
for $f(x)$.
The requirement that $alg_{ls}(x)$ compute an answer as accurately as though 
it were backward stable raises its bit complexity only by a factor
at most quadratic in $\chi(n)$, i.e. polynomial in $\log n$.
\end{lemma}

\begin{proof}
A backward stable algorithm for $f(x)$ running with
machine precision $\varepsilon_{bs}$ would have 
relative error bound $O( \varepsilon_{bs} ) \kappa_f(x) = \tau$. 
A relative error bound is only meaningful when it is less than 1,
so we may assume $\tau < 1$. Taking logarithms yields
the number of bits $b_{bs}$ of precision needed:
\begin{equation}
b_{bs} = \log_2 \frac{1}{\varepsilon_{bs}} = \log_2 \frac{1}{\tau} + \log_2 \kappa_f(x) 
+ O(1) \; \; .
\end{equation}
Recall that each arithmetic operation costs at most $O( b_{bs}^2)$ bit operations
and as few as $O( b_{bs} \log b_{bs} \log \log b_{bs} )$ if fast techniques are
used \cite{SchonhageStrassen}. 

To make the actual error bound for $alg_{ls}(x)$ as small as $\tau$
means we have to choose $\varepsilon_{ls}$
to satisfy $O( \varepsilon_{ls} ) \kappa_f^{\chi(n)}(x) = \tau$.
Again taking logarithms yields the number of bits $b_{ls}$ 
of precision needed:
\begin{equation}
b_{ls} = \log_2 \frac{1}{\varepsilon_{ls}} 
= \log_2 \frac{1}{\tau} + \chi(n) \cdot \log_2 \kappa_f(x) + O(1)
\leq \chi(n) b_{bs} + O(1)
\end{equation}
This raises the cost of each arithmetic operation in $alg_{ls}(x)$ 
by a factor of at most $O(\chi^2(n))$ as claimed.

Thus, if $alg_{ls}(x)$ were
backward stable and performed $O(n^c)$ arithmetic operations,
it would cost at most $O(n^c b_{bs}^2)$ bit operations to get
a relative error $\tau < 1$.
Logarithmic stability raises its cost to 
at most $O(n^c \chi^2(n) b_{bs}^2)$ bit operations to get
the same relative error.
\end{proof}

\subsection{Recursive Triangular Matrix Inversion}
\label{sec_RTriInv}

First we apply these ideas to triangular matrix inversion, based on the formula
\[
T^{-1} = \bmat{cc} T_{11} & T_{12} \\ 0 & T_{22} \emat^{-1}
   = \bmat{cc} T_{11}^{-1} & -T_{11}^{-1} \cdot T_{12} \cdot T_{22}^{-1} \\
                0 & T_{22}^{-1} \emat
\]
where $T_{11}$ and $T_{22}$ are $\frac{n}{2}$-by-$\frac{n}{2}$ and
inverted using the same formula recursively.
The cost of this well-known algorithm \cite{BorodinMunro75,Heller78}
is $$cost(n) = 2 \, cost(n/2) + 2MM(n/2,n/2,n/2) = O(n^{\gamma}).$$
Its error analysis in \cite{higham95parallel} (Method A in Section 6)
used the stronger componentwise bound \cite{higham96}[eqn. (3.13)]
that holds for conventional matrix multiplication 
(as opposed to (\ref{gen_bound}))
but nevertheless concluded that the method was
not as stable as the conventional method.
(The motivation for considering this algorithm in \cite{higham95parallel}
was not fast matrix multiplication but parallelism, which also leads
to many block algorithms.) 

To show this algorithm is logarithmically stable,
we do a first order error analysis for the absolute error $err(T^{-1},n)$
in the computed inverse of the $n$-by-$n$ matrix $T$.
We use the fact that in computing the product of two matrices $C = A \cdot B$
that have inherited errors $err(A,n)$ and $err(B,n)$ from prior computations,
we may write
\begin{eqnarray}
\label{eqn_MulBnd}
err(C,n) & = & \mu(n) \varepsilon \|A\| \cdot \|B\| 
               \hspace*{.58in} {\rm \ from\ matrix\ multiplication} \nonumber \\
         &   & + \|A\| \cdot err(B,n)
               \hspace*{.5in} {\rm \ amplifying\ the\ error\ in\ }B{\rm \ by\ }\|A\| \\
         &   & + err(A,n) \cdot \|B\|
               \hspace*{.5in} {\rm \ amplifying\ the\ error\ in\ }A{\rm \ by\ }\|B\| \nonumber 
\end{eqnarray}
We will also use the facts that 
$\|T_{ii}\| \leq \|T\|$ 
(and $\|T_{ii}^{-1}\| \leq \|T^{-1}\|$) 
since $T_{ii}$ is a submatrix of $T$ 
(and $T_{ii}^{-1}$ is a submatrix of $T^{-1}$).
Therefore the condition number
$\kappa (T_{ii}) \equiv \|T_{ii}\| \cdot \|T_{ii}^{-1}\| \leq \kappa (T)$.
Now let $err(n')$ be a bound for the normwise error in the inverse
of any $n'$-by-$n'$ diagonal subblock of $T$ encountered during the algorithm.
Applying (\ref{eqn_MulBnd}) systematically to the recursive algorithm yields the
following recurrences bounding the growth of $err(n)$.
(Note that we arbitrarily decide to premultiply $T_{12}$ by $T_{11}^{-1}$ first.)
\begin{eqnarray*}
err(T_{ii}^{-1},n/2)              & \leq &  err(n/2) 
             \hspace{3cm}               ... {\rm \ from\ inverting\ }T_{11} {\rm \ and\ }T_{22} \\
err(T_{11}^{-1} \cdot T_{12},n/2) & \leq &  \mu(n/2) \varepsilon \|T_{11}^{-1}\| \cdot \|T_{12}\|
                                            + err(T_{11}^{-1},n/2) \|T_{12}\| \\
        & & \hspace{4.3cm}  ... {\rm \ from\ multiplying\ }T_{11}^{-1} \cdot T_{12} \\
                                  & \leq & \mu(n/2) \varepsilon \|T^{-1}\| \cdot \|T\| + err(n/2) \|T\| \\
err((T_{11}^{-1} \cdot T_{12}) \cdot T_{22}^{-1},n/2) & \leq &  
                  \mu(n/2) \varepsilon \|T_{11}^{-1} \cdot T_{12}\| \cdot \|T_{22}^{-1}\| \\
       &      & + err(T_{11}^{-1} \cdot T_{12},n/2) \cdot \|T_{22}^{-1}\| \\
       &      & + \|T_{11}^{-1} \cdot T_{12}\| \cdot err(T_{22}^{-1},n/2) \\
       &      &     \hspace{4.3cm} ... {\rm \ from\ multiplying\ }(T_{11}^{-1} \cdot T_{12}) \cdot T_{22}^{-1} \\
       & \leq & \mu(n/2) \varepsilon \|T^{-1}\| \cdot \|T\| \cdot \|T^{-1}\| \\
        &  &    + (\mu(n/2) \varepsilon \|T^{-1}\| \cdot \|T\| + err(n/2) \|T\|) \cdot \|T^{-1}\| \\
        &  &    + \|T^{-1}\| \cdot \|T\| \cdot err(n/2) \\
err(T^{-1},n) & \leq & err(T_{11}^{-1},n/2) + err(T_{22}^{-1},n/2) + 
                       err((T_{11}^{-1} \cdot T_{12}) \cdot T_{22}^{-1},n/2) \\ 
              &  \leq &  2 err(n/2) \\
          & &   + \mu(n/2) \varepsilon \|T^{-1}\| \cdot \|T\| \cdot \|T^{-1}\| \\
          & &   + (\mu(n/2) \varepsilon \|T^{-1}\| \cdot \|T\| + err(n/2) \|T\|) \cdot \|T^{-1}\|\\
          & &   + \|T^{-1}\| \cdot \|T\| \cdot err(n/2) \\
       & \leq & 2(\kappa(T) +1) err(n/2) + 
                2\mu(n/2) \varepsilon \kappa(T) \|T^{-1}\| 
\end{eqnarray*}
Solving the resulting recurrence for $err(n)$ \cite{CLR}[Thm. 4.1]
yields
\begin{eqnarray*}
err(n) & = & 2(\kappa(T) +1) err(n/2) + 
                2\mu(n/2) \varepsilon \kappa(T) \|T^{-1}\| \\
        & = & 
        O(\mu(n/2) \varepsilon \kappa(T) (2(\kappa(T)+1))^{\log_2 n} \|T^{-1}\|)
\end{eqnarray*}
showing that the algorithm is logarithmically stable as claimed.
\ignore{
This formula tells us how small $\varepsilon$ has to be to attain a desired absolute
error bound $err(n) \leq \varepsilon_0 \|T^{-1}\|$. Computing $log_2 \frac{1}{\varepsilon}$
tells us how many bits of precision we need to carry, and so the complexity of
each operation:
\[
\log_2 \frac{1}{\varepsilon} = O(\log_2 \frac{1}{\varepsilon_0}  
                               + (\log_2 n) \cdot (\log \kappa (T)) )
\]
In comparison, a conventional algorithm with an absolute error bound
$O(\varepsilon \kappa(T) \|T^{-1}\|)$ would need
$\log_2 \frac{1}{\varepsilon} = O(\log_2 \frac{1}{\varepsilon_0} + \log \kappa (T) )$
bits, only a factor $O(\log_2 n)$ lower.
Thus, the extra complexity from using high precision to get
an answer as accurate as the conventional algorithm only contributes
a polylog($n$) factor to the overall complexity of triangular matrix inversion
as claimed.
} 

\subsection{Recursive Dense Matrix Inversion}

A similar analysis may be applied to inversion of symmetric positive definite
matrices using an analogous well-known divide-and-conquer formula:
\begin{eqnarray*}
H & = & \bmat{cc} A & B \\ B^T & C \emat =
\bmat{cc} I & 0 \\ B^TA^{-1} & I \emat \cdot 
\bmat{cc} A & B \\ 0 & S \emat 
\hspace{.5in} {\rm where\ } S = C - B^T A^{-1} B \\
\Longrightarrow
H^{-1} & = & 
\bmat{cc} A^{-1} & -A^{-1}BS^{-1} \\ 0 & S^{-1} \emat \cdot
\bmat{cc} I & 0 \\ -B^TA^{-1} & I \emat  =
\bmat{cc} A^{-1} + A^{-1}B S^{-1} B^T A^{-1} & -A^{-1}BS^{-1} \\ 
                    -S^{-1}B^T A^{-1} & S^{-1} \emat 
\end{eqnarray*}
To proceed, we need to state the algorithm more carefully,
also deriving a recurrence for the cost $C(n)$:

\begin{tabbing}
jnkk \= jnkk \= jnkk \= jnkkkkkkkkkkkkkkkkkkkkkkkk \= jnk \= jnk \kill
\> function $Hi = RecursiveInv(H,n)$ ... invert $n$-by-$n$ s.p.d. matrix 
$H$ recursively \\
\> \> if (n=1) then \\
\> \> \> $Hi = 1/H$ \\
\> \> else \\
\> \> \> $Ai = RecursiveInv(A,n/2)$  \>... cost = $C(n/2)$ \\
\> \> \> $AiB = Ai \cdot B$   \>... cost = $MM(n/2)$ \\
\> \> \> $BAiB = B^T \cdot AiB$   \>... cost = $MM(n/2)$ \\
\> \> \> $S = C - BAiB$ \>... cost = $(n/2)^2$ \\
\> \> \> $Si = RecursiveInv(S,n/2)$ \>... cost = $C(n/2)$ \\
\> \> \> $AiBSi = AiB \cdot Si$   \>... cost = $MM(n/2)$ \\
\> \> \> $AiBSiBAi = AiBSi \cdot (AiB)^T$ \>... cost = $MM(n/2)$ \\
\> \> \> $Hi_{11} = Ai + AiBSiBAi$ \>... cost = $(n/2)^2$ \\
\> \> \> return $Hi = [[Hi_{11}, -AiBSi];[ (-AiBSi)^T, Si]]$ \\
\> \> endif
\end{tabbing}

Assuming $MM(n) = O(n^{\gamma})$ for some $2 < \gamma \leq 3$, it is easy to
see that the solution of the cost recurrence
$C(n) = 2C(n/2) + 4MM(n/2) + n^2/2 = O(MM(n))$ as desired.

For the rest of this section the matrix norm $\| \cdot \|$ will
denote the 2-norm (maximum singular value).
To analyze the error we exploit the Cauchy Interlace Theorem which 
implies first that the eigenvalues of $A$ interlace the eigenvalues of $H$, 
so $A$ can be no worse conditioned than $H$, and second that
the eigenvalues of $S^{-1}$ (and so of $S$) interlace the eigenvalues
of $H^{-1}$ (and so of $H$, resp.), so $S$ can also be no worse
conditioned than $H$. Letting $\lambda$ and $\Lambda$ denote
the smallest and largest eigenvalues of $H$, resp., we also get
that $\| B^TA^{-1}B \| \leq \|C\| + \|S\| \leq 2 \Lambda$
and $\|A^{-1} BS^{-1}B^T A^{-1} \| \leq \|A^{-1}\| + \|H^{-1}\| \leq
2/ \lambda$, all of which inequalities we will need below.

As before, we use the induction hypothesis that
$err(n')$ bounds the error in the inverse
of any $n'$-by-$n'$ diagonal block computed by the
algorithm (including the errors in computing the
block, if it is a Schur complement, as well as inversion).
In particular we assume $err(A^{-1},n/2) \leq err(n/2)$. 
Then we get
\begin{align*}
err(AiB,n/2) & \leq 
     \mu(n/2) \cdot \varepsilon \cdot \frac{1}{\lambda} \cdot \Lambda
     + err(n/2) \cdot \Lambda  
 \hspace{3.3cm}  ... {\rm using\ }(\ref{eqn_MulBnd}){\rm \ with\ no\ error\ in\ }B \\
err(BAiB,n/2) & \leq  \mu(n/2) \cdot \varepsilon \cdot \Lambda \cdot
                       \frac{\Lambda}{\lambda} +
                       \Lambda \cdot err(AiB,n/2) 
\hspace{2.4cm}  ... {\rm using\ }(\ref{eqn_MulBnd}){\rm \ with\ no\ error\ in\ }B & \\
              & \leq  2 \mu(n/2) \cdot \varepsilon \cdot \frac{\Lambda^2}{\lambda} 
                + \Lambda^2 \cdot err(n/2) \\
err(S,n/2) & \leq  \sqrt{\frac{n}{2}} \varepsilon \cdot 
             \Lambda + err(BAiB,n/2) \\
           & \approx  err(BAiB,n/2) 
\\ 
err(Si,n/2) & \leq  err(n/2) + \frac{1}{\lambda^2} \cdot err(S,n/2) 
\hspace{2.5cm}  ... {\rm using\ }          (S + \delta S)^{-1} \approx S^{-1} - S^{-1} \delta S S^{-1} & \\
           & \leq  err(n/2) + 2 \mu(n/2) \cdot \varepsilon \cdot \frac{\Lambda^2}{\lambda^3} 
                + \frac{\Lambda^2}{\lambda^2} \cdot err(n/2) 
\end{align*}
\begin{align*}
err(AiBSi,n/2) & \leq  \mu(n/2) \cdot \varepsilon \cdot 
                        \frac{\Lambda}{\lambda} \cdot \frac{1}{\lambda} +
                        \frac{\Lambda}{\lambda} \cdot err(Si,n/2) +
                        err(AiB,n/2)  \cdot \frac{1}{\lambda} 
\hspace{2.4cm}  ... {\rm using\ }(\ref{eqn_MulBnd}) & \\
                  & \leq  \mu(n/2) \cdot \varepsilon \cdot 
                           (\frac{2 \Lambda}{\lambda^2} + \frac{2 \Lambda^3}{\lambda^4}) +
                           (\frac{2 \Lambda}{\lambda} + \frac{\Lambda^3}{\lambda^3}) err(n/2)\\
err(AiBSiBAi,n/2) & \leq  \mu(n/2) \cdot \varepsilon \cdot
                           \frac{1}{\lambda} \cdot \frac{\Lambda}{\lambda} +
                           \frac{1}{\lambda} \cdot err(AiB,n/2) +
                           err(AiBSi,n/2) \cdot \frac{\Lambda}{\lambda} 
\hspace{1.6cm}  ... {\rm using\ }(\ref{eqn_MulBnd}) \\
                  & \leq  \mu(n/2) \cdot \varepsilon \cdot 
                           (\frac{2 \Lambda}{\lambda^2} + \frac{2 \Lambda^2}{\lambda^3} +
                            \frac{2 \Lambda^4}{\lambda^5}) +
                           (\frac{\Lambda}{\lambda} + \frac{2 \Lambda^2}{\lambda^2} + 
                            \frac{\Lambda^4}{\lambda^4}) err(n/2) \\
err(Hi_{11},n/2)  & \leq  \sqrt{\frac{n}{2}} \varepsilon 
                    \cdot \frac{1}{\lambda} + err(A^{-1},n/2) + 
                           err(AiBSiBAi,n/2) \\
                  & \approx  \mu(n/2) \cdot \varepsilon \cdot 
                           (\frac{2 \Lambda}{\lambda^2} + \frac{2 \Lambda^2}{\lambda^3} +
                            \frac{2 \Lambda^4}{\lambda^5}) +
                           (1 + \frac{\Lambda}{\lambda} + \frac{2 \Lambda^2}{\lambda^2} + 
                            \frac{\Lambda^4}{\lambda^4}) err(n/2) \\
err(Hi,n)          & \leq  err(Hi_{11},n/2) + err(AiBSi,n/2) + err(Si,n/2) \\
                  & \leq  \mu(n/2) \cdot \varepsilon \cdot
                           (4 \frac{\Lambda}{\lambda^2} +
                            4 \frac{\Lambda^2}{\lambda^3} +
                            2 \frac{\Lambda^3}{\lambda^4} +
                            2 \frac{\Lambda^4}{\lambda^5}) +
                            err(n/2) \cdot
                           (2 + 3 \frac{\Lambda}{\lambda} +
                            3 \frac{\Lambda^2}{\lambda^2} +
                              \frac{\Lambda^3}{\lambda^3} +
                             \frac{\Lambda^4}{\lambda^4})
\end{align*}


This yields a recurrence for the error, where we write $\kappa = \frac{\Lambda}{\lambda}$:
\begin{eqnarray*}
  err(n) & \leq & \mu(n/2) \cdot \varepsilon \cdot
                           (4 \frac{\Lambda}{\lambda^2} +
                            4 \frac{\Lambda^2}{\lambda^3} +
                            2 \frac{\Lambda^3}{\lambda^4} +
                            2 \frac{\Lambda^4}{\lambda^5}) +
                            err(n/2) \cdot
                           (2 + 3 \frac{\Lambda}{\lambda} +
                            3 \frac{\Lambda^2}{\lambda^2} +
                              \frac{\Lambda^3}{\lambda^3} +
                             \frac{\Lambda^4}{\lambda^4}) \\
         & \leq & 12 \cdot \mu(n/2) \cdot \varepsilon  \cdot \frac{\kappa^4}{\lambda}
                  + 10 \cdot \kappa^4 \cdot err(n/2) 
\end{eqnarray*}
Solving this recurrence, we get
\begin{equation}
\label{eqn_InvLogStable}
err(n) = O( \varepsilon \mu(n) \kappa^4 (10 \kappa^4)^{\log_2 n} \lambda^{-1} ) 
       = O( \varepsilon \mu(n) n^{\log_2 10} \kappa^{4 + 4 \log_2 n} \|H^{-1}\| )
\end{equation}
showing that recursive inversion of a symmetric positive definite matrix is
logarithmically stable.

\ignore{
An analogous but more involved first order error analysis for $err(H^{-1},n)$ 
leads to the recurrence (we omit details)
\begin{eqnarray*}
err(n) & \leq & (9 \kappa^4(H) + 4 \kappa^2(H)) err(n/2) + 16 \mu(n/2) \varepsilon
                  \kappa^4(H) \|H^{-1} \| \\
       & \leq & (13 \kappa^4(H)) err(n/2) + 16 \mu(n/2) \varepsilon \kappa^4(H) \|H^{-1}\| \\
       &  =   & O( \mu (n/2) n^{\log_2 13} \varepsilon \|H^{-1}\| \kappa^{4\log_2(n)}(H)
\end{eqnarray*}
This again shows that the number of bits to get an accurate answer grows proportionally to 
$O(\log_2 \frac{1}{\varepsilon_0} + (\log_2 n) \cdot (\log \kappa (T)) )$,
only $O(\log_2 n)$ times as many as a conventional algorithm.
} 

To invert a general matrix we may use $A^{-1} = A^T \cdot (A \cdot A^{T})^{-1}$.
Forming $A \cdot A^T$ only squares $A$'s condition number, and 
first order error analysis shows the errors contributed from the 
two matrix multiplications can only increase the exponent of
$\kappa$ in (\ref{eqn_InvLogStable}) by doubling it and adding a 
small constant.
Thus we may also draw the conclusion 
that general matrix inversion is logarithmically stable.
The same reasoning applies to solving $Ax=b$ by multiplying 
$x = A^T \cdot (A \cdot A^{T})^{-1} \cdot b$.

Finally, we return to our claim in the introduction:

\begin{theorem}
\label{thm_1}
If we can multiply $n$-by-$n$ matrices in 
$O(n^{\omega + \eta})$ arithmetic operations 
then we can invert matrices stably
in $O(n^{\omega + \eta})$ bit operations.
Conversely, if we can invert matrices 
stably in $O(n^{\omega + \eta})$ bit operations
(resp. exactly in $O(n^{\omega + \eta})$ arithmetic operations)
then we can multiply matrices
stably in $O(n^{\omega + \eta})$ bit operations
(resp. exactly in $O(n^{\omega + \eta})$ arithmetic operations).
\end{theorem}

\begin{proof}
We have just proven the first claim, where we rely on logarithmic
stability of inversion to bound the number of bit operations.
\ignore{
We have just proved the forward implication: 
matrix multiplication in $O(n^{\omega + \eta})$ {\em bit} operations 
certainly implies
matrix multiplication in $O(n^{\omega + \eta})$ {\em arithmetic} operations 
(since each arithmetic operation uses multiple bit operations).
This in turn implies stable (in the sense of (\ref{gen_bound}))
matrix multiplication in $O(n^{\omega + \eta})$ {\em arithmetic} operations 
by \cite{FastMatMulIsStable}. This in turn implies we can
invert matrices stably (in the sense of Lemma~\ref{lemma_LogStable})
in $O(n^{\omega + \eta})$ {\em bit} operations.
}

For the converse implications, we simply use
\[
\bmat{ccc} I & A & 0 \\ &  I & B \\ & & I \emat^{-1} =
\bmat{ccc} I & -A & A \cdot B \\ &  I & -B \\ & & I \emat \\.
\]
Clearly, inverting the matrix on the left 
exactly in $O(n^{\omega + \eta})$ arithmetic operations
lets us extract the product $A \cdot B$.
Given only a logarithmically stable inversion routine,
we can make the condition number near 1 by
scaling $A$ and $B$ to have norms near 1,
implying that the above block matrices are very well
conditioned, and the inverse can be computed accurately
without extra precision.
\end{proof}

It is tempting to summarize this theorem by saying
``matrix multiplication is possible in
$O(n^{\omega + \eta})$ operations if and only if stable inversion is,''
but the difference between counting bit operations and arithmetic
operations requires a more careful statement
(a bound on the number of arithmetic operations can be used
to bound the number of bit operations, but not conversely,
since bit operations may conceivably not organize themselves
into easily recognized arithmetic operations).

\section{Simultaneous Speed and Backward Stability of QR and LU}
\label{sec_Simultaneous}

We show that QR decomposition can be implemented stably
and as fast as matrix multiplication. We exploit the fact
that linear equation solving and determinant computation
as well as solving least squares problems
can be reduced to QR decomposition to make the same statements
about these linear algebra operations. Similar statements
can be made about LU decomposition, under slightly stronger
assumptions about pivot growth than the conventional algorithm.

\subsection{Fast and Stable QR Decomposition}
\label{sec_QRStable}

We now describe in more detail the following recursive variation of
a conventional QR algorithm \cite{BischofVanLoan87},
which was presented in \cite{ElmrothGustavson2000}.
Let $A$ be an $n$-by-$m$ matrix with $n \geq m$.
The function $[R, W, Y] = QRR(A,n,m)$ will return an $m$-by-$m$ upper triangular
matrix $R$, an $n$-by-$m$ matrix $W$, and an $m$-by-$n$
matrix $Y$ with the following properties:
(1) $Q^T = I - W Y$ is an $n$-by-$n$ orthogonal matrix, 
(2) each column of $W$ has unit 2-norm, 
(3) each row of $Y$ has 2-norm equal to 2, 
(4) $A = Q \cdot [R;zeros(n-m,m)]$ is the QR decomposition of $A$
(here and later we use MATLAB notation).

\begin{tabbing}
jnk \= jnk \= jnk \= jnk \= jnk \kill
\> function $[R, W, Y] = QRR(A)$ ... $A$ is $n$-by-$m$, with $n \geq m$ \\
\> \> if ($m=1$) then \\
\> \> \> compute $W$ and $Y$ in the conventional way as a Householder transformation
         \cite[sec. 19.1]{higham96}, \\
\> \> \> \> with the normalization that $\|W\|_2 = 1$, $\|Y\|_2 = 2$ and 
           $R = \pm \|A\|_2$ \\
\> \> else \\
(a) \> \> \> $[R_L, W_L, Y_L] = QRR(A(1:n, 1: \mtf ))$ \\
             \> \> \> \> ... {\rm \ compute\ QR\ decomposition\ of\ left\ half\ of\ }$A$ \\
(b) \> \> \> $A(1:n, \mtf + 1:m) = A(1:n, \mtf +1:m) -
          W_L \cdot (Y_L \cdot A(1:n, \mtf +1:m))$ \\
             \> \> \> \> ... {\rm \ multiply\ right\ half\ of\ }$A${\rm \ by\ }$Q^T$ \\
(c) \> \> \> $[R_R, W_R, Y_R] = QRR(A(\mtf +1:n, \mtf +1:m ))$ \\
             \> \> \> \> ... {\rm \ compute\ QR\ decomposition\ of\ right\ half\ of\ }$A$ \\
(d) \> \> \> $X = W_L - [zeros( \mtf , \mtf );W_R \cdot (Y_R \cdot W_L( \mtf +1:n,1: \mtf )]$ \\
               \> \> \> \> ... {\rm \ multiply\ two\ }$Q${\rm \ factors} \\
\> \> \> $R = [[R_L, A(1:\mtf, \mtf+1:m)];[zeros(\mtc , \mtf ), R_R]]$ \\
\> \> \> $W = [X, [zeros(\mtf , \mtc ); W_R]]$ \\
\> \> \> $Y = [Y_L ; [zeros(\mtc , \mtf ), Y_R]]$ \\
\> \> endif \\
\end{tabbing}

The proof of correctness is induction based on the identity
$(I - [0;W_R][0,Y_R]) \cdot (I - W_L Y_L) = I - WY$
as in Section~2 above.
For the complexity analysis we assume $m$ is a power of 2:
\begin{eqnarray*}
cost(n,m) & = & cost(n,\mt ) 
                \hspace{3.46in} ... {\rm cost\ of\ line\ (a)} \\
          &   &  + MM( \mt , n, \mt ) + MM( n, \mt , \mt ) +
                 n \mt
                \hspace{1.585in} ... {\rm cost\ of\ line\ (b)} \\
        &   &  + cost(n - \mt , \mt )
                \hspace{3.03in} ... {\rm cost\ of\ line\ (c)} \\
          &   &  + MM( \mt, n- \mt , \mt ) + MM( n- \mt , \mt , \mt ) 
                  + (n- \mt ) \mt
                \hspace{.5in} ... {\rm cost\ of\ line\ (d)} 
\end{eqnarray*}
\begin{eqnarray*}
 &\leq & 2 cost(n, \mt ) + 8 \frac{n}{m} MM( \mt , \mt , \mt ) 
                  + O(nm) \\
          &\leq & 2 cost(n, \mt ) + O(nm^{\gamma -1}) \\
          & =   & O(nm^{\gamma -1})
                \hspace{3.46in} ... {\rm assuming\ } \gamma > 2
\end{eqnarray*}
When $n=m$, this means the complexity is $O(n^{\gamma})$ as desired.


This algorithm submits to an analogous backward error analysis
as in \cite{BischofVanLoan87} or \cite{higham96}[sec. 19.5], which
we sketch here for completeness. 

\begin{lemma}
\label{lemma_QRStability}
The output of $[R,W,Y]=QRR(A)$ satisfies
$(I - WY + \delta Q^T)(A + \delta A) = [R;zeros(n-m,m)]$ where
$Q^T \equiv I - WY + \delta Q^T$ satisfies $QQ^T = I$ exactly,
$\| \delta Q^T \| = O(\varepsilon)$, 
and $\| \delta A \| = O( \varepsilon ) \|A\|$.
(Here we let $O()$ absorb all factors of the form $n^c$.)
\end{lemma}

\begin{proof}
We use proof by induction. 
The base case ($m=1$) may be
found in \cite{higham96}[sec. 19.3].
Let $A_L = A(1:n,1:\mtf)$ and $A_R = A(1:n,\mtf+1:m)$.
From the induction hypothesis applied to step (a) of QRR we
have 
\[
(I - W_L Y_L + \delta Q_L^T)( A_L + \delta A_L ) 
= [R_L;zeros(n- \mtf ,\mtf )]
\; \; \; {\rm with} \; \; \;
Q_L^T \equiv I - W_L Y_L + \delta Q_L^T \; \; ,
\]
$Q_L^T Q_L = I$,
$\| \delta Q_L^T \| = O( \varepsilon )$ and 
$\| \delta A_L \| = O( \varepsilon ) \|A\|$.
Application of error bound (\ref{gen_bound}) to step (b) yields
\[
A_{R,new} = A_R - W_L \cdot (Y_L \cdot A_R) + \delta A_{R,1}
= Q_L^T (A_R + \delta \hat{A}_{R,1}) 
\; \; \; {\rm with} \; \; \;
\delta \hat{A}_{R,1} = -Q_L \cdot \delta Q_L^T \cdot A_R + Q_L \cdot \delta A_{R,1}
\]
so $\| \delta \hat{A}_{R,1} \| = O( \varepsilon ) \|A\|$.
Write $A_{R,new} = [A_{R,1}; A_{R,2}]$ where $A_{R,1}$ is $\mtf$-by-$\mtc$.
The induction hypothesis applied to step (c) yields
\[
(I - W_R Y_R + \delta Q_R^T)(A_{R,2} + \delta A_{R,2}) = [R_R;zeros(n-m,\mtc]
\; \; \; {\rm with} \;  \; \;
Q_R^T \equiv I - W_R Y_R + \delta Q_R^T
\]
$Q_R^T Q_R = I$,
$\| \delta Q_R^T \| = O( \varepsilon )$, and
$\| \delta A_{R,2} \| = O( \varepsilon ) \|A\|$.
Combining expressions we get
\[
\bmat{cc} I & 0 \\ 0 & Q_R^T \emat \cdot Q_L^T \cdot (A + \delta A)
 = \bmat{cc} R_L & A_{R,1} \\ 0 & R_R \\ 0 & 0 \emat
\]
where
\[
\delta A = \left[ \delta A_L , \delta \hat{A}_{R,1} + 
   Q_L \cdot \bmat{c} zeros(\mtf , \mtc ) \\ \delta A_{R,2} \emat \right]
\]
satisfies $\| \delta A \| = O( \varepsilon ) \|A\|$.
Finally, repeated application of bound~(\ref{gen_bound}) to
step (d) shows that 
$X = X_{true} + \delta X$ with $\| \delta X \| = O( \varepsilon )$,
$W = W_{true} + [\delta X, zeros(n,\mtc )]$, and
\begin{eqnarray*}
Q^T & \equiv & \bmat{cc} I & 0 \\ 0 & Q_R^T \emat \cdot Q_L^T \\
& = & 
       \bmat{cc} I & 0 \\ 0 & I - W_R Y_R + \delta Q_R^T \emat 
      \cdot
      (I - W_L Y_L + \delta Q_L^T) \\
& = & I - W_{true} Y + \delta \hat{Q}^T \\
& = & I - WY + \delta Q^T
\end{eqnarray*}
with $QQ^T = I$,
$\| \delta \hat{Q}^T \| = O( \varepsilon )$ and
$\| \delta Q^T \| = O( \varepsilon )$ as desired.
\end{proof}


Armed with an $A=QR$ decomposition, we can easily solve the linear system
$Ax = b$ stably via $x = R^{-1}Q^T b$ straightforwardly 
in another $O(n^2)$ operations, or solve a least squares problem stably.
Furthermore det$(A) = (-1)^n \prod_i R_{ii}$ is also easily computed.
In summary, high speed and numerical stability are achievable
simultaneously.

\subsection{Fast and Stable LU Decomposition}
\label{sec_LUStable}

There is an analogous algorithm for LU decomposition \cite{Toledo97}.
However, in order to update the right half of the matrix after doing
the LU decomposition of the left half, it appears necessary to
invert a lower triangular matrix,
namely the upper left corner of the $L$ factor, whose inverse
is then multiplied by the upper right corner of $A$ to get the
upper right corner of $U$.
As described in the last section,
triangular matrix inversion seems to be only logarithmically stable.
However, because of pivoting, one is guaranteed that $L_{ii} = 1$
and $|L_{ij}| \leq 1$, so that $\kappa(L)$ is generally small.
Thus as long as $L$ is sufficiently well conditioned 
then LU decomposition can also be done
stably and as fast as matrix multiplication. Now we sketch the details,
omitting the implementation of pivoting, since it does not contribute
to the complexity analysis:

\begin{tabbing}
jnk \= jnk \= jnk \= jnk \= jnk \kill
\> function $[L, U] = LUR(A)$ ... $A$ is $n$-by-$m$, with $n \geq m$ \\
\> \> if (m=1) then \\
\> \> \> $L = A/A(1)$, $U = A(1)$ \\
\> \> else \\
(a) \> \> \> $[L_L, U_L] = LUR(A(1:n, 1: \mtf ))$ \\
             \> \> \> \> ... {\rm \ compute\ LU\ decomposition\ of\ left\ half\ of\ }$A$ \\
(b) \> \> \> $A(1:\mtf , \mtf + 1:m) = (L_L(1:\mtf , 1:\mtf ))^{-1} \cdot A(1:\mtf , \mtf +1:m)$; \\
             \> \> \> \> ... {\rm \ update\ upper\ right\ corner\ of\ }$A$ \\
(c) \> \> \> $A(\mtf +1:n , \mtf + 1:m) = A(\mtf +1:n , \mtf +1:m) - $ \\
    \> \> \> \> $L_L(\mtf +1:n , 1:\mtf ) \cdot A(1:\mtf , \mtf + 1:m)$; \\
             \> \> \> \> ... {\rm \ update\ Schur\ complement} \\
(d) \> \> \> $[L_R, U_R] = LUR(A(\mtf +1:n, \mtf +1:m ))$ \\
             \> \> \> \> ... {\rm \ compute\ LU\ decomposition\ of\ right\ half\ of\ }$A$ \\
(e) \> \> \> $L = [L_L, [zeros( \mtf , \mtc ) ; L_R ] ]$; \\
(f) \> \> \> $U = [[U_L, A(1:\mtf , \mtf +1 : m)];[zeros(\mtc , \mtf ) , U_R ] ]$; \\
\> \> endif \\
\end{tabbing}

For the complexity analysis we assume $m$ is a power of 2 as before:

\begin{eqnarray*}
cost(n,m) & = & cost(n,\mt ) 
                \hspace{1in} ... {\rm cost\ of\ line\ (a)} \\
          &   &  + O(MM( \mt )) 
                \hspace{1in} ... {\rm cost\ of\ line\ (b)} \\
        &   &  + MM(n - \mt , \mt , \mt ) + (n- \mt ) \mt 
                \hspace{1in} ... {\rm cost\ of\ line\ (c)} \\
          &   &  + cost(n-\mt , \mt )
                \hspace{1in} ... {\rm cost\ of\ line\ (d)} \\
 &\leq & 2 cost(n, \mt ) + 2 \frac{n}{m} MM( \mt , \mt , \mt ) 
                  + O(nm + MM( \mt )) \\
          &\leq & 2 cost(n, \mt ) + O(nm^{\gamma -1}) \\
          & =   & O(nm^{\gamma -1})
                \hspace{1in} ... {\rm assuming\ } \gamma > 2
\end{eqnarray*}
When $n=m$, this means the complexity is $O(n^{\gamma})$ as desired.

Now we establish backward stability under the assumption that
$L$ (and so every diagonal block of $L$) is sufficiently well
conditioned (and its norm sufficiently close to 1) 
that the error in the computed matrix from step (b)
is bounded in norm by $O( \varepsilon \|A\|)$:

\begin{lemma}
\label{lemma_LUStability}
If $L$ in the output of $[L,U]=LUR(A)$ is
sufficiently well conditioned, then
$L \cdot U = A + \delta A$ where $\| \delta A \| = O( \varepsilon ) \|A\|$.
(Here we let $O()$ absorb all factors depending on
$\| L \|$, $\| L^{-1} \|$, and $n$.
We also assume without loss of generality that the 
rows of $A$ are in the correct pivot order.)
\end{lemma}

\begin{proof}
We use proof by induction. 
The base case ($m=1$) is straightforward.
Let 
$A_L = A(1:n , 1:\mtf)$,
$L_{L,1} = L(1:\mtf , 1:\mtf)$,
$L_{L,2} = L(\mtf + 1:n , 1:\mtf)$,
$A_{R,1} = A(1:\mtf , \mtf + 1 : m)$,
and
$A_{R,2} = A(\mtf + 1:n , \mtf + 1 : m)$.
Then from the induction hypothesis applied to step (a),
$L_L \cdot U_L = A_L + \delta A_L$ with 
$\| \delta A_L \| = O( \varepsilon \|A\| )$.
From step (b) and the assumptions about $L$,
the updated value of $A_{R,1}$ is given by
\[
A'_{R,1} = L_{L,1}^{-1} \cdot A_{R,1} + \delta A'_{R,1}
\; \; \; {\rm with} \; \; \;
\| \delta A'_{R,1} \| = O( \varepsilon \|A\| ) \; \; .
\]
From step (c) the updated value of $A_{R,2}$ is given by
\[
A'_{R,2} = A_{R,2} - L_{L,2} \cdot A'_{R,1} + \delta A'_{R,2}
\; \; \; {\rm with} \; \; \; 
\| \delta A'_{R,2} \| = O( \varepsilon \|A\| ) \; \; .
\]
From the induction hypothesis applied to step (d) we get 
\[
L_R \cdot U_R = A'_{R,2} + \delta A''_{R,2}
\; \; \; {\rm with} \; \; \;
\| \delta A''_{R,2} \| = O( \varepsilon \|A\| ) \; \; .
\]
Combining these results yields
\[
L \cdot U = A + \delta A
\; \; \; {\rm with} \; \; \; 
\delta A = [\delta A_L , [L_{L,1} \delta A'_{R,1} ; \delta A''_{R,2} + \delta A'_{R,2} ] ] \; \; ,
\]
so $\| \delta A \| = O( \varepsilon \|A\| )$ as desired.
\end{proof}

The assumption that $L$ is sufficiently well-conditioned is a variation
on the usual assumption that pivot growth is limited, since pivot
growth is bounded by $\|L^{-1}\|$ (with the norm depending on how
pivot growth is measured), and $\|L\|_1$ is at most $n$.

\subsection{Columnwise Backward Error}

The error analysis of conventional $O(n^3)$ algorithms for the QR and LU
decomposition actually yield somewhat stronger results than normwise
backward stability: they are normwise backward stable {\em column-by-column}.
This means, for example, that LU is the exact factorization of $A + \delta A$
where the $i$-th column of $\delta A$ is small in norm compared to the 
$i$-th column of $A$. As stated, our algorithm does not have this property,
since the fast matrix multiplication algorithm can and probably will
``smear'' errors between columns. But there is a simple fix to avoid this,
and get the same columnwise error bound as the standard algorithms:
(1) Preprocess $A$ by dividing each column by its norm (say the infinity norm);
save the values of the norms for step (3).
(2) Compute the fast QR (or LU) factorization of the scaled $A$.
(3) Multiply the $i$th-column of $R$ (or of $U$) by the norm of the $i$-th
column of $A$.

It is easy to see that this additional work costs only $O(n^2)$, and makes
the backward error in column $i$ proportional to the norm of column $i$.

More generally, one could improve bound~(\ref{gen_bound}) either to
$|C_{comp,ij} - C_{ij}| \leq \mu(n) \epsilon \|A(i,:)\| \, \|B(:,j)\|$
or to
$\|C_{comp} - C\| \leq \mu(n) \epsilon \| \, |A| \cdot |B| \, \|$
by appropriately scaling rows and/or columns of $A$ and $B$ before
multiplying them, and unscaling $C_{comp}$ afterwards if necessary.

\section{Fast and Stable Randomized Rank Revealing URV}
\label{sec_URV}

We also show how to implement a rank revealing URV decomposition
based on QR decomposition stably and fast; 
this will be required for solving eigenvalue problems in the next section. 
Our rank revealing algorithm will be {\em randomized}, i.e. it will work
with high probability. 
As we will see in the next section, this is adequate for our eigenvalue algorithm.

Given a (nearly) rank deficient matrix $A$, our goal is to
quickly and stably compute a factorization $A = URV$ where
$U$ and $V$ are orthogonal and $R$ is upper triangular,
with the property that it {\em reveals the rank} in the
following sense: Let $\sigma_1 \geq \cdots \geq \sigma_n$ be
the singular values of $A$. Then
(1) with high probability $\sigma_{\min}(R(1:r,1:r))$ is a good approximation of  
$\sigma_r$, and (2) assuming
there is a gap in the singular values ($\sigma_{r+1} \ll \sigma_r$)
and that $R(1:r,1:r)$ is not too ill-conditioned, then
with high probability
$\sigma_{\max}(R(r+1:n,r+1:n))$ is a good approximation 
of $\sigma_{r+1}$.
This is analogous to other definitions of rank-revealing decompositions
in the literature \cite{chan87,hongpan,chandipsen,bischof96,gueisenstat2,stewart93}, 
with the exception of its randomized nature. 

The algorithm is quite simple
(RURV may be read ``randomized URV'')

\begin{tabbing}
jnk \= jnk \= jnk \= jnk \= jnk \kill
\> function $[U, R, V] = RURV(A)$ ... $A$ is $n$-by-$n$ \\
\> \> generate a random matrix $B$ whose entries are independent, identically distributed \\
\> \> \> Gaussian random variables with mean 0 and standard deviation 1 (i.i.d. N(0,1)) \\
(a) \> \> $[V,R] = QRR(B)$ \hspace{.5in} ... $V$ is a random orthogonal matrix \\
(b) \> \> $\hat{A} = A \cdot V^T$ \\
(c) \> \> $[U,R]=QRR(\hat{A})$ 
\end{tabbing}

\ignore{
Thus $U \cdot R = \hat{A} = A \cdot V^T$, so $U \cdot R \cdot V = A$.
The cost of RURV is one matrix multiplication and two calls to QRR,
plus $O(n^2)$ to form $B$, so $O(n^{\omega + \eta})$ altogether.
The matrix $V$ is a random orthogonal matrix in the following sense...
{(\em Ioana, please fill in)}.
It remains to prove that this is a rank-revealing decomposition
with high probability
(for simplicity we restrict ourselves to real matrices, 
although the analysis easily extends to complex matrices):

\begin{theorem}
In exact arithmetic, the $R$ matrix produced by $RURV(A)$ satisfies 
two conditions. 
First, 
\[
f \cdot \sigma_r \leq \sigma_{\min} (R(1:r,1:r)) \leq \sqrt{\sigma_r^2 + \sigma_{r+1}^2}
\leq \sqrt{2} \cdot \sigma_r
\]
where $f$ is a random variable 
equal to the smallest singular value of an $r$-by-$r$ submatrix of
a random $n$-by-$n$ orthogonal matrix, 
and so has distribution XXX {\em (Ioana?)}.
Second, assuming $\sigma_{r+1} < f \sigma_r$,
\[
\sigma_{\max} (R(r+1:n,r+1:n))
\leq 
3 \sigma_{r+1} \cdot
\frac{f^{-4} \cdot \left( \frac{\sigma_1}{\sigma_r} \right)^3}{1 - \frac{\sigma_{r+1}^2}{f^2 \sigma_r^2}}
\]
This means that if there is a large gap in the singular values ($\sigma_{r+1} \ll \sigma_r$)
and $R(1:r,1:r)$ is not too ill-conditioned ($\sigma_r$ is not $\ll \sigma_1$),
then with high probability the output of RURV is rank-revealing.
\end{theorem}

} 

Thus $U \cdot R = \hat{A} = A \cdot V^T$, so $U \cdot R \cdot V = A$.
The cost of RURV is one matrix multiplication and two calls to QRR,
plus $O(n^2)$ to form $B$, so $O(n^{\omega + \eta})$ altogether.
The matrix $V$ has \emph{Haar distribution} 
\cite{muirhead82a},
i.e. it is distributed uniformly over the set of $n$-by-$n$ orthogonal matrices.
For information on efficient generation of such matrices, see
\cite{anderson87a,stewart80a}.

It remains to prove that this is a rank-revealing decomposition
with high probability
(for simplicity we restrict ourselves to real matrices, 
although the analysis easily extends to complex matrices):

\begin{lemma} \label{prob_bound}
Let $f$ be a random variable equal to the smallest singular value of an $r$-by-$r$ submatrix of a Haar distributed random $n$-by-$n$ orthogonal matrix. Assume that $r$ is ``large'' (i.e, grows to $\infty$ as some function of $n$; no assumptions are made on the growth speed). Let $a>0$ be a positive constant. Then there is a constant $c>1$ such that, as soon as $r$ and $n$ are large enough,
\[
\mbox{Pr} \left [ f < \frac{1}{r^{a+1} \sqrt{n}} \right ] ~ \leq ~ \frac{c}{r^{a}}~.
\]
\end{lemma}

\begin{proof}
Recall that the $n$-by-$n$ Haar distribution has the following important 
property: any column, as well as the transpose of any row, is 
uniformly distributed over the $n-1$ unit sphere. As such, 
without loss of generality, we can restrict ourselves to the study of the 
leading $r$-by-$r$ submatrix.

Denote by $V$ the orthogonal matrix, and let $U = V(1:r, 1:r)$. We will 
assume that $V$ came from the QR factorization of a $n$-by-$n$ matrix $B$ of 
i.i.d. Gaussians (for simplicity, as in $RURV$), 
and thus $V = B R^{-1}$. Moreover, $U = 
B(1:r, 1:r) (R(1:r, 1:r))^{-1}$. Therefore 
\begin{eqnarray*} 
f := \sigma_{\mbox{min}}(U) & \geq & \sigma_{\mbox{min}}(B(1:r, 1:r)) 
\cdot 
\sigma_{\mbox{min}}((R(1:r, 1:r))^{-1})~~, \\ 
& \geq &  \frac{\sigma_{\mbox{min}}(B(1:r, 1:r))}{\sigma_{\mbox{max}}(R(1:r, 1:r))}~~, \\ & \geq & \frac{\sigma_{\mbox{min}}(B(1:r, 1:r))}{ 
\sigma_{\mbox{max}}(B(1:n, 1:r))}~~,\\ & \geq & 
\frac{\sigma_{\mbox{min}}(B(1:r, 1:r))}{||B(1:n, 1:r)||_{F}}~~, 
\end{eqnarray*} where $||~||_F$ denotes the Frobenius norm.

Thus we shall have that
\[
\mbox{Pr} \left [ f < \frac{1}{r^{a+1} \sqrt{n} } \right ] \leq \mbox{Pr} \left [ 
\frac{\sigma_{\mbox{min}}(B(1:r, 1:r))}{||B(1:n, 1:r)||_{F}} <  \frac{1}{r^{a+1} \sqrt{n} } \right 
]~~.
\]

We now use the following bound:
\begin{eqnarray} \label{bound}
\mbox{Pr} \left [
\frac{\sigma_{\mbox{min}}(B(1:r, 1:r))}{||B(1:n, 1:r)||_{F}} < \frac{1}{r^{a+1 } \sqrt{n} } \right
] \leq \mbox{Pr} \left [ \sigma_{\mbox{min}}(B(1:r, 1:r)) < \frac{2}{r^{a+1/2}} \right ] +  \mbox{Pr} \left [ ||B(1:n, 1:r)||_{F} > 2 \sqrt{rn} \right ] 
~.
\end{eqnarray}

The limiting (asymptotical) distribution (as $r \rightarrow \infty$) of $r \cdot \sigma_{\mbox{min}}(B(1:r, 1:r))^2$ has been computed in Corollary~3.1 of 
\cite{edelman88} and shown to be given by 
\[
f(x) = \frac{1 + \sqrt{x}}{2 \sqrt{x}} ~e^{-(x/2+\sqrt{x})}~.
\]
The convergence was shown to be very fast; in particular there exists a (small) constant $c_0$ such that 
\[
\mbox{Pr} \left[  r^2 \sigma_{\mbox{min}}(B(1:r, 1:r)) < x \right] \leq c_0 \int_0^x f(t)dt < c_0 \sqrt{x} ~,
\]
for any $x>0$. After the appropriate change of variables, it follows that there is a constant $c_1$ such that  
\begin{eqnarray} \label{unu}
\mbox{Pr} \left [ \sigma_{\mbox{min}}(B(1:r, 1:r)) < \frac{2}{r^{a+1/2}} \right ] \leq \frac{c_1}{r^{a}}~,
\end{eqnarray}
for all $r$.

On the other hand, the distribution of the variable $||B(1:n, 1:r)||_{F}$ is $\chi_{nr}$, with $\chi$ being the square root of the $\chi^2$ variable. As the probability density function for $\chi_{rn}$ is
\[
g_{rn}(x) = \frac{1}{2^{rn/2-1} \Gamma\left( \frac{rn}{2} \right)} ~~x^{\frac{rn}{2}-1}~e^{-x^2/2} ~,
\]
and simple calculus gives the bound
\begin{eqnarray} \label{doi}
 \mbox{Pr} \left [ ||B(1:n, 1:r)||_{F} > 2 \sqrt{rn} \right ] \leq e^{-\sqrt{rn}/2}~,
\end{eqnarray}
for all $r$ and $n$.

From \eqref{bound}, \eqref{unu}, and \eqref{doi} we obtain the statement of the lemma.
\end{proof}

Lemma~\ref{prob_bound} implies that, as long as $r$ grows with $n$,
the chance that $f$ is small is itself small, certainly less than
half, which is all we need for a randomized algorithm to work in
a few trials with high probability.

\begin{theorem}
In exact arithmetic, the $R$ matrix produced by $RURV(A)$ satisfies  
the following two conditions. First,
\[
f \cdot \sigma_r \leq \sigma_{\min} (R(1:r,1:r)) \leq \sqrt{\sigma_r^2 + \sigma_{r+1}^2}
\leq \sqrt{2} \cdot \sigma_r~,
\]
where $f$ is a random variable 
equal to the smallest singular value of an $r$-by-$r$ submatrix of
a random $n$-by-$n$ orthogonal matrix. Second, assuming $\sigma_{r+1} < f \sigma_r$,
\[
\sigma_{r+1} \leq
\sigma_{\max} (R(r+1:n,r+1:n))
\leq 
3 \sigma_{r+1} \cdot
\frac{f^{-4} \cdot \left( \frac{\sigma_1}{\sigma_r} \right)^3}{1 - \frac{\sigma_{r+1}^2}{f^2 \sigma_r^2}}
\]
Given Lemma \ref{prob_bound}, which says it is unlikely for $f$ to be small,
this means that if there is a large gap in the singular values of 
$A$ ($\sigma_{r+1} \ll \sigma_r$) and $R(1:r,1:r)$ is not too ill-conditioned 
($\sigma_r$ is not $\ll \sigma_1$),
then with high probability the output of RURV is rank-revealing.
\end{theorem}

\begin{proof}
In this proof $\| Z \|$ will denote the largest singular value of $Z$.
Let $A = P \cdot \Sigma \cdot Q^T = P \cdot \diag ( \Sigma_1 , \Sigma_2 ) \cdot Q^T$ 
be the singular value decomposition of $A$, where
$\Sigma_1 = \diag ( \sigma_1 , ... , \sigma_r )$ and
$\Sigma_2 = \diag ( \sigma_{r+1} , ... , \sigma_n )$.
Let $V^T$ be the random orthogonal matrix in the RURV algorithm.
Note that $X \equiv Q^T \cdot V^T$ has the same probability distribution as $V^T$.
The intuition is simply that a randomly chosen subspace (spanned
by the leading $r$ columns of $V^T$) is unlikely to contain
vectors nearly orthogonal to another $r$ dimensional subspace
(spanned by the leading $r$ columns of $Q$),
i.e. that the leading $r$-by-$r$ submatrix of $X$ is unlikely 
to be very ill-conditioned.
Write 
$X = [ X_1 , X_2 ]$
where $X_1 = \bmat{c} X_{11} \\ X_{21} \emat$ has $r$ columns,
$X_2 = \bmat{c} X_{12} \\ X_{22} \emat$ has $n-r$ columns,
$X_{11}$ and $X_{12}$ have $r$ rows, and
$X_{21}$ and $X_{22}$ have $n-r$ rows.
Then 
\[
\sigma_{\min}(R(1:r,1:r)) 
 = \sigma_{\min} \left( \bmat{cc} \Sigma_1 \cdot X_{11} \\ \Sigma_2 \cdot X_{21} \emat \right)
 \geq \sigma_{\min} (\Sigma_1 \cdot X_{11})
 \geq \sigma_r \cdot \sigma_{\min} (X_{11})
 = \sigma_r \cdot f \; \; .
\]
where $f \equiv \sigma_{\min} (X_{11})$ is a random variable with
distribution described in Lemma~\ref{prob_bound}.
Also
\[
\sigma^2_{\min} (R(1:r,1:r)) 
  = \sigma_{\min} (X_{11}^T \Sigma_1^2 X_{11} + X_{21}^T \Sigma_2^2 X_{21} )
  \leq \sigma_{\min} (X_{11}^T \Sigma_1^2 X_{11} ) + \sigma_{\max} (X_{21}^T \Sigma_2^2 X_{21} )
  \leq \sigma_r^2 + \sigma_{r+1}^2 \;\; .
\]
Now let $\Sigma \cdot X = [ \Sigma \cdot X_1 , \Sigma \cdot X_2 ] \equiv [ Y_1 , Y_2 ]$.
Then the nonzero singular values of $R(r+1:n,r+1:n)$ are identical to singular values
of the projection of $Y_2$ on the orthogonal complement of the column space
of $Y_1$, namely of the matrix $C = (I - Y_1 (Y_1^T Y_1)^{-1} Y_1^T) Y_2$.
Write 
\[
(Y_1^T Y_1)^{-1} = (X_{11}^T \Sigma_1^2 X_{11} + X_{21}^T \Sigma_2^2 X_{21})^{-1}
 \equiv (S - \delta S)^{-1}
 = S^{-1} + \sum_{i=1}^{\infty} S^{-1} ( \delta S \cdot S^{-1})^i
\]
assuming the series converges.
Assuming the product of $\sigma^2_{r+1} \geq \| \delta S \|$ and
of $\frac{1}{f^2 \sigma_r^2} \geq \|S^{-1}\|$ is less than 1, we have
\[
\hat{Y} \equiv (Y_1^T Y_1)^{-1} = (X_{11}^T \Sigma_1^2 X_{11})^{-1} + E
 = X_{11}^{-1} \Sigma_1^{-2} X_{11}^{-T} + E
\]
where 
$\|E\| \leq (\frac{\sigma_{r+1}^2}{f^4 \sigma_r^4}) / (1 - \frac{\sigma_{r+1}^2}{f^2 \sigma_r^2})$
and $\| \hat{Y} \| \leq \frac{1}{f^2 \sigma_r^2} + \|E\| 
                   \leq \frac{1}{f^2 \sigma_r^2 - \sigma_{r+1}^2}$.
Now we compute
\begin{eqnarray*}
C & = & (I - Y_1 (Y_1^T Y_1)^{-1} Y_1^T) Y_2 \\
& = &
\bmat{cc} I - \Sigma_1 X_{11} \hat{Y} X_{11}^T \Sigma_1 & 
          -\Sigma_1 X_{11} \hat{Y} X_{21}^T \Sigma_2    \\
          -\Sigma_2 X_{21} \hat{Y} X_{11}^T \Sigma_1    &
          I - \Sigma_2 X_{21} \hat{Y} X_{21}^T \Sigma_2 \emat \cdot
\bmat{c} \Sigma_1 X_{12} \\ \Sigma_2 X_{22} \emat \\
& = &
\bmat{cc} -\Sigma_1 X_{11} E  X_{11}^T \Sigma_1         & 
          -\Sigma_1 X_{11} \hat{Y} X_{21}^T \Sigma_2    \\
          -\Sigma_2 X_{21} \hat{Y} X_{11}^T \Sigma_1    &
          I - \Sigma_2 X_{21} \hat{Y} X_{21}^T \Sigma_2 \emat \cdot
\bmat{c} \Sigma_1 X_{12} \\ \Sigma_2 X_{22} \emat \\
& = &
\bmat{c} -\Sigma_1 X_{11} E  X_{11}^T \Sigma_1^2 X_{12} \\
         -\Sigma_2 X_{21} \hat{Y} X_{11}^T \Sigma_1^2 X_{12} \emat
 +
\bmat{c}
          -\Sigma_1 X_{11} \hat{Y} X_{21}^T \Sigma_2    \\
          I - \Sigma_2 X_{21} \hat{Y} X_{21}^T \Sigma_2 \emat \cdot
          \Sigma_2 X_{22} \\
& \equiv &
    \bmat{c} Z_1 \\ Z_2 \emat + Z_3 \cdot \Sigma_2 X_{22}
\end{eqnarray*}
Thus $\|C\| \leq \| Z_1 \| + \| Z_2 \| + \| Z_3 \| \cdot \| \Sigma_2 \| \cdot \| X_{22} \|
\leq \| Z_1 \| + \| Z_2 \| + \sigma_{r+1}$.
Next, 
\[
\| Z_1 \| 
\leq 
\| \Sigma_1 \| \cdot \| X_{11} \| \cdot \| E \| \cdot
\| X_{11} \| \cdot \| \Sigma_1^2 \| \cdot \| X_{12} \| 
\leq 
\sigma_1^3 \| E \| 
\leq 
\frac{\sigma_1^3 \frac{\sigma_{r+1}^2}{f^4 \sigma_r^4}}{1 - \frac{\sigma_{r+1}^2}{f^2 \sigma_r^2}}
\leq 
\sigma_{r+1} \cdot
\frac{f^{-4} \cdot \left( \frac{\sigma_1}{\sigma_r} \right)^3}{1 - \frac{\sigma_{r+1}^2}{f^2 \sigma_r^2}}
\]
and
\[
\| Z_2 \| 
\leq 
\| \Sigma_2 \| \cdot \| X_{21} \| \cdot \| \hat{Y} \| \cdot
\| X_{11} \| \cdot \| \Sigma_1^2 \| \cdot \| X_{12} \| 
\leq 
\sigma_1^2 \sigma_{r+1} \| \hat{Y} \|
\leq
\frac{ \sigma_1^2 \sigma_{r+1} }{f^2 \sigma_r^2 - \sigma_{r+1}^2}
\]
which is smaller than the bound on $\| Z_1 \|$, as is $\sigma_{r+1}$.
Altogether, we then get
\[
\| C \| \leq 
3 \sigma_{r+1} \cdot
\frac{f^{-4} \cdot \left( \frac{\sigma_1}{\sigma_r} \right)^3}{1 - \frac{\sigma_{r+1}^2}{f^2 \sigma_r^2}}
\]
as desired.

\end{proof}

\begin{corollary}
Suppose $A$ has rank $r<n$. Then (in exact arithmetic) 
the $r$ leading columns of $U$
from $[U,R,V] = RURV(A)$ span the column space of $A$ with
probability 1.
\end{corollary}

\begin{lemma}
\label{lemma_URVStability}
In the presence of roundoff error, the computed
output $[U,R,V]=RURV(A)$ satisfies 
\newline
$A + \delta A = \hat{U} \cdot R \cdot \hat{V}$
where 
$\hat{U}$ and $\hat{V}$ are exactly orthogonal matrices,
$\| \delta A \| = O( \varepsilon ) \|A\|$,
$\|U - \hat{U}\| = O( \varepsilon )$, and
$\|V - \hat{V}\| = O( \varepsilon )$.
\end{lemma}

\begin{proof}
Applying Lemma~\ref{lemma_QRStability} to step (a) yields 
$V = \hat{V} - \delta V$ where $\hat{V} \cdot \hat{V}^T = I$
and $\| \delta V \| = O( \varepsilon )$.
Applying error bound~(\ref{gen_bound}) to step (b) yields
$\hat{A} = A \cdot V^T + \delta A_1$ where
$\| \delta A_1 \| = O( \varepsilon \|A\| )$.
Applying Lemma~\ref{lemma_QRStability} to step (c) yields 
$\hat{U} \cdot R = \hat{A} + \delta A_2$ where
$\hat{U} \cdot \hat{U}^T = I$, $\delta U = \hat{U} - U$
satisfies $\| \delta U \| = O( \varepsilon )$, and
$\| \delta A_2 \| = O( \varepsilon \|A\| )$.
Combining these identities yields
\[
\hat{U} \cdot R = (A - A \cdot \delta V^T \cdot \hat{V} +
(\delta A_1 + \delta A_2) \cdot \hat{V} ) \cdot \hat{V}^T
\equiv (A + \delta A) \cdot \hat{V}^T
\]
where $\| \delta A \| = O( \varepsilon \|A\|)$ as desired.
\end{proof}

Lemma~\ref{lemma_URVStability} shows that $RURV(A)$ computes
a rank revealing factorization of a matrix close to $A$,
which is what is needed in practice (see the next section).
(We note that merely randomizing the order of the columns
of $A$ is not good enough: consider the case of an $n$-by-$n$
matrix of rank $r$ where $n-r+1$ columns are all multiples
of one another.)
The question remains of how to recognize success, that a
rank-revealing factorization has in fact been computed.
(The same question arises for conventional rank-revealing
QR, with column pivoting, which can fail, rarely, on 
matrices like the Kahan matrix.)
This will be done as part of the eigenvalue algorithm.


\section{Eigenvalue Problems} 
\label{sec_Eigen}

To show how to solve eigenvalue problems quickly and stably,
we use an algorithm from \cite{baidemmelgu94}, modified
slightly to use only the randomized rank revealing decomposition 
from the last section.  
As described in Section~\ref{sec_Schur},
it can compute either an invariant subspace of a 
matrix $A$, or a pair of left and right deflating 
subspaces of a regular matrix pencil $A - \lambda B$, using
only QRR, RURV and matrix multiplication. 
Applying it recursively, and with some assumptions
about partitioning the spectrum, we can compute a (generalized)
Schur form stably in $O(n^{\omega + \eta})$ arithmetic operations.
Section~\ref{sec_SymEigSVD} discusses the special case
of symmetric matrices and the singular value decomposition,
where the previous algorithm is enough to stably compute eigenvectors
(or singular vectors) as well in
$O(n^{\omega + \eta})$ operations.
But to compute eigenvectors of a nonsymmetric matrix (or pencil)
in $O(n^{\omega + \eta})$ operations is more difficult: 
Section~\ref{sec_Sylvester} gives a logarithmically stable
algorithm $SylR$ to solve the (generalized) Sylvester equation,
which Section~\ref{sec_Eigenvectors} in turn uses in 
algorithm $EVecR$ for eigenvectors. 
However, the $EVecR$ is only logarithmically stable in a
weak sense: the accuracy of any computed eigenvector
may depend on the condition numbers of other, worse
conditioned eigenvectors. We currently see no way
to compute each eigenvector with an error proportional
to its own condition number other than by the conventional
$O(n^2)$ algorithm (involving the solution of a triangular
system of equations), for a cost of $O(n^3)$ to compute
all the eigenvectors this accurately.

\subsection{Computing (generalized) Schur form} 
\label{sec_Schur}

The stability and convergence analysis of the algorithm is subtle, 
and we refer to \cite{baidemmelgu94} for details. (Indeed, not
even the conventional Hessenberg QR algorithm has a convergence
proof \cite{day96}.) Our goal here is to show that the algorithm
in \cite{baidemmelgu94} can be implemented as stably as described
there, and in $O(n^{\omega + \eta})$ operations. The only
change required in the algorithm is replacing use of the
QR decomposition with pivoting, and how rank is determined: 
instead we use the RURV decomposition, and determine rank
by a direct assessment of backward stability described below.
Since this only works with high probability,
we will have to loop, repeating the decomposition with a
different random matrix, until a natural stopping criterion is 
met. The number of iterations will be low since the
chance of success at each step is high.

Rather than explain the algorithm in 
\cite{baidemmelgu94} in detail, we briefly describe a similar, simpler
algorithm in order to motivate how one might
use building blocks like matrix multiplication, inversion,
and QR decomposition to solve the eigenvalue problem.
This simpler algorithm is based on the {\em matrix sign-function}
\cite{roberts}:
Suppose $A$ has no eigenvalues with zero imaginary part,
and let $A = S \cdot \diag (J_+, J_- ) \cdot S^{-1}$ be
its Jordan Canonical form, where $J_+$ consists
of the Jordan blocks for the $r$ eigenvalues with positive real part,
and $J_-$ for the negative real part.
Then define $\sign (A) \equiv S \cdot \diag (+I_r, -I_{n-r}) \cdot S^{-1}$.
One can easily see that $P_+ = \frac{1}{2} (\sign(A) + I ) = 
S \cdot \diag (+I_r, 0) \cdot S^{-1}$ is the
spectral projector onto the invariant subspace ${\cal S}_+$
of $A$ for $J_+$.  
Now perform a rank revealing QR factorization of $P_+$, yielding
an orthogonal matrix $Q$ whose leading $r$ columns span ${\cal S}_+$.
Therefore 
\[
Q^T A Q = \bmat{cc} A_{11} & A_{12} \\ 0 & A_{22} \emat
\]
is a {\em block Schur factorization}:
it is orthogonally similar to $A$, and the $r$-by-$r$ matrix
$A_{11}$ has all the eigenvalues with positive imaginary part,
and 
$A_{22}$ has all the eigenvalues with negative imaginary part.

We still need a method to compute $\sign (A)$. Consider
Newton's method applied to find the zeros of $f(x) = x^2 - 1$,
namely $x_{i+1} = \frac{1}{2} (x_i + x_i^{-1})$.
Since the signs of the real parts of $x_i$ and $x_i^{-1}$, and
so of $x_{i+1}$, are identical, Newton can only converge to $\sign (\Re x_0)$.
Global and eventual quadratic convergence follow by using
the Cayley transform to change variables to $\hat{x}_i = \frac{x_i-1}{x_i+1}$:
This both maps the left (resp. right) half plane to the 
exterior (resp. interior) of the unit circle, and Newton 
to $\hat{x}_{i+1} = \hat{x}_i^2$, whose convergence is apparent.
We therefore use the same iteration for matrices:
$A_{i+1} = \frac{1}{2}(A_i + A_i^{-1})$
(see \cite{roberts}).
One can indeed show this converges globally and
ultimately quadratically to $\sign (A_0)$.
This reduces computing the sign function to matrix
inversion.

To compute a more complete Schur factorization,
we must be able to apply this algorithm recursively
to $A_{11}$ and $A_{22}$, and so divide their
spectra elsewhere than along the imaginary axis.
By computing a Moebius transformation
$\hat{A} = (\alpha A + \beta I) \cdot (\gamma A + \delta I)^{-1}$,
we can transform the imaginary axis to an arbitrary
line or circle in the complex plane. So by computing
a rank-revealing QR decomposition of $\frac{1}{2}( \sign ( \hat{A} ) +I )$,
we can compute an invariant subspace for the eigenvalues inside or
outside any circle, or in any halfspace. To automate the choice
of these regions, one may use 2-dimensional bisection, or quadtrees,
starting with a rectangle guaranteed to contain all eigenvalues (using
Gershgorin bounds), and repeatedly dividing it into smaller subrectangles,
stopping division when the rectangle contains too few eigenvalues
or has sufficiently small perimeter.

The same ideas may be extended to the generalized eigenvalue problem,
computing left and right deflating subspaces of the regular matrix pencil
$A - \lambda B$.

Backward stability (though not progress) may be guaranteed by
checking whether the computed $A_{21}$ in
$Q^T A Q = \bmat{cc} A_{11} & A_{12} \\ A_{21} & A_{22} \emat$
is sufficiently small in norm. See \cite{baidemmel92a,baidemmel97,HussEtAl,PRISM}
for a more complete description of eigenvalue algorithms
based on the matrix sign function.

The algorithm in \cite{baidemmelgu94},
which in turn is based on earlier algorithms
of Bulgakov, Godunov and Malyshev
\cite{godunov86,bulgakov88,malyshev89,malyshev92,malyshev93},
avoids all matrix inverses, either to compute a basis
for an invariant (or deflating) subspace, or to 
split the spectrum along an arbitrary line or circle in the 
complex plane. It only uses QR decompositions and matrix
multiplication to compute a matrix whose columns span 
the desired invariant subspace, followed by a rank-revealing
QR decomposition to compute an orthogonal matrix $Q$ whose
leading columns span the subspace, and an upper triangular
matrix $R$ whose rank determines the dimension of the subspace.
Stability is enforced as above by checking to make
sure $A_{21}$ is sufficiently small in norm, rejecting the
decomposition if it is not. 

To modify this algorithm to use our RURV decomposition,
we need to have a loop to repeat the RURV decomposition
until it succeeds (which takes few iterations with 
high probability), and determine the rank
(number of columns of $A_{21}$).
Rather than determine the rank from $R$, we
choose it to minimize the norm of $A_{21}$,
by computing $\| A_{21} \|_S \equiv \sum_{ij} |A_{21,ij}|$
for all possible choices of rank, in just $O(n^2)$ operations:
\begin{tabbing}
jnk \= jnk \= jnk \= jnk \= jnk \kill
\> repeat \\
\> \> compute $Q$ from RURV \\
\> \> compute $\hat{A} = Q^T A Q$ \\
\> \> for $i=1$ to $n-1$ \\
\> \> \> ColSum($i$) $= \sum_{j=i+1}^n |\hat{A}_{ji}|$ \\
\> \> \> Rowsum($i+1$) $= \sum_{j=1}^{i} |\hat{A}_{i+1,j}|$ \\
\> \> endfor \\
\> \> NormA21($1$) = ColSum($1$) \\
\> \> for $i=2$ to $n-1$ \\
\> \> \> NormA21($i$) = NormA21($i-1$) $+$ ColSum($i$) $-$ RowSum($i$) \\
\> \> \> \> ... NormA21($i$) $= \| \hat{A}(i+1:n,1:i) \|_S$ \\
\> \> endfor \\
\> \> Let $r$ be the index of the minimum value of NormA21($1:n-1$) \\
\> until NormA21($r$) small enough, or too many iterations
\end{tabbing}

We need the clause ``too many iterations'' in the algorithm to
prevent infinite loops, and to account for the possibility that
the algorithm was unable to split the spectrum at all. This could
be because all the eigenvalues had positive imaginary part (or
were otherwise on the same side of the dividing line or circle),
or were too close to the dividing line or circle that the algorithm used.
The error analysis of the algorithm in \cite{baidemmelgu94}
is otherwise identical.

Finally we need to confirm that the overall complexity of the
algorithm, including updating the Schur form and accumulating
all the orthogonal transformations to get a complete
Schur form $A = QTQ^T$, costs just $O(n^{\omega + \eta})$.
To this end, suppose the original matrix $A$ is $n$-by-$n$, and the 
dimension of a diagonal submatrix $\bar{A}$ encountered during 
divide-and-conquer is $\bar{n}$. Then computing an $\bar{n}$-by-$\bar{n}$
$\bar{Q}$ to divide $\bar{A}$ once takes $O(\bar{n}^{\omega + \eta})$
operations as described above, 
applying $\bar{Q}$ to the rest of $A$ costs at most
$2 \cdot MM(\bar{n} , \bar{n} , n ) = O(n \bar{n}^{\omega + \eta - 1})$
operations, and accumulating $\bar{Q}$ into the overall $Q$ costs
another $MM(\bar{n} , \bar{n} , n ) = O(n \bar{n}^{\omega + \eta - 1})$.
Letting $cost(\bar{n})$ be the cost of all work associated with $\bar{A}$
then yields the recurrence 
\[
cost(\bar{n}) 
= 2 cost( \frac{\bar{n}}{2} ) + O(\bar{n}^{\omega + \eta}) + O(n \bar{n}^{\omega + \eta - 1})
= 2 cost( \frac{\bar{n}}{2} ) + O(n \bar{n}^{\omega + \eta - 1})
\]
whose solution is $cost(\bar{n}) = O(n \bar{n}^{\omega + \eta - 1})$
or
$cost(n) = O(n^{\omega + \eta})$ as desired.

\subsection{Symmetric Matrices and the SVD}
\label{sec_SymEigSVD}

When the matrix is symmetric, or is simply known to have
real eigenvalues, simpler alternatives to the matrix sign
function are known that only involve matrix multiplication
\cite{bischofledermansuntsao,hussledermantsaozhang93,PRISM},
and to which the techniques described here may be applied.
Symmetry can also be enforced stably in the algorithm described
above by replacing each computed matrix $Q^TAQ$ by its
symmetric part. The bisection technique described above
to locate eigenvalues of general matrices obviously simplifies 
to 1-dimensional bisection of the real axis.

The SVD of $A$ can be reduced to the symmetric eigenproblem 
either (1) for $\bmat{cc} 0 & A \\ A^T & 0 \emat$ or (2) for
$AA^T$ and $A^TA$. But in either case, the notion of backward stability 
for the computed singular vectors would have to be modified slightly to 
account for possible difficulties in computing them.
We can avoid this difficulty, and get a fully backward stable SVD algorithm,
by separately computing orthogonal $Q_L$ and $Q_R$ whose leading 
$r$ columns (nearly) span a left (resp. right) singular subspace of $A$, 
and forming $Q_L^T A Q_R$. This should be (nearly) block diagonal,
letting us continue with divide-and-conquer.
$Q_L$ (resp. $Q_R$) would be computed by applying our earlier
algorithm to compute an orthogonal matrix whose leading $r$ columns
span an eigenspace of $AA^T$ (resp. $A^TA$, for the same subset of the
spectrum). Stability despite squaring $A$ requires double precision,
a cost hidden by the big-O analysis.
The algorithm would also check
to see if $Q_L^T A Q_R$ is close enough to block diagonal, for any 
block size $r$, to enforce stability.

\subsection{Solving the (generalized) Sylvester Equation}
\label{sec_Sylvester}

To compute an invariant subspace of 
$\bmat{cc} A & C \\ 0 & B \emat$
for the eigenvalues of $B$ and
spanned by the columns of $\bmat{c} R \\ I \emat$ 
we are lead to the equation 
\begin{equation}
\label{eqn_InvariantSS}
\bmat{cc} A & C \\ 0 & B \emat
\cdot
\bmat{c} R \\ I \emat
=
\bmat{c} R \\ I \emat
\cdot B
\end{equation}
or the {\em Sylvester equation}
$A R - R B = -C$ to solve for $R$. When $A$ and $B$ 
are upper triangular as in Schur form,
this is really a permuted triangular system of equations
for the entries of $R$, where the diagonal
entries of the triangular matrix are all 
possible differences $A_{ii} - B_{jj}$,
so the system is nonsingular precisely when
the eigenvalues of $A$ and $B$ are distinct.

Similarly, to compute a right (resp. left) deflating subspace
of $\bmat{cc} A & C \\ 0 & B \emat - \lambda 
\bmat{cc} \bar{A} & \bar{C} \\ 0 & \bar{B} \emat$
for the eigenvalues of $B - \lambda \bar{B}$ and
spanned by the columns of $\bmat{c} R \\ I \emat$ 
(resp. $\bmat{c} L \\ I \emat$) 
we are led to the equations
\[
\bmat{cc} A & C \\ 0 & B \emat \cdot \bmat{c} R \\ I \emat
= \bmat{c} L \\ I \emat \cdot B
\; \; \; {\rm and} \; \; \; 
\bmat{cc} \bar{A} & \bar{C} \\ 0 & \bar{B} \emat \cdot \bmat{c} R \\ I \emat
= \bmat{c} L \\ I \emat \cdot \bar{B}
\]
or the {\em generalized Sylvester equation}
$A R - L B = -C$, $\bar{A} R - L \bar{B} = - \bar{C}$ to solve for $R$ and $L$. 

To derive a divide-and-conquer algorithm for the Sylvester equation,
write it as
\begin{equation}
\label{eqn_Syl}
\bmat{cc} A_{11} & A_{12} \\ 0      & A_{22} \emat \cdot
\bmat{cc} R_{11} & R_{12} \\ R_{21} & R_{22} \emat -
\bmat{cc} R_{11} & R_{12} \\ R_{21} & R_{22} \emat \cdot
\bmat{cc} B_{11} & B_{12} \\ 0      & B_{22} \emat =
- \bmat{cc} C_{11} & C_{12} \\ C_{21} & C_{22} \emat 
\end{equation}
where all submatrices are dimensioned conformally.
Multiplying this out we get the four equations
\begin{eqnarray}
A_{22} R_{21} - R_{21} B_{11} & = & - C_{21} 
     \label{eqn_Syl1} \\
A_{11} R_{11} - R_{11} B_{11} & = & - C_{11} - A_{12} R_{21} 
     \label{eqn_Syl2} \\
A_{22} R_{22} - R_{22} B_{22} & = & - C_{22} + R_{21} B_{12} 
     \label{eqn_Syl3} \\
A_{11} R_{12} - R_{12} B_{22} & = & - C_{12} + R_{11} B_{12} - A_{12} R_{22}
    \label{eqn_Syl4} 
\end{eqnarray}
We recognize this as four smaller Sylvester equations, where
(\ref{eqn_Syl1}) needs to be solved first for $R_{21}$,
which lets us evaluate the right hand sides of 
(\ref{eqn_Syl2}) and (\ref{eqn_Syl3}), and finally  (\ref{eqn_Syl4}).

This idea is captured in the following algorithm, where for
simplicity we assume all matrices are square and of the
same dimension, a power of 2:

\begin{tabbing}
jnk \= jnk \= jnk \= jnk \= jnk \kill
\> function $R = SylR(A,B,C)$ ... all matrices are $n$-by-$n$ \\
\> \> if ($n=1$) then \\
\> \> \> $R = -C/(A-B)$ \\
\> \> else \hspace*{3.5in} ... use notation from equation (\ref{eqn_Syl}) \\
(a) \> \> \> $R_{21} = SylR(A_{22},B_{11},C_{21})$
    \hspace*{1.9in} ... solve (\ref{eqn_Syl1}) \\
(b) \> \> \> $R_{11} = SylR(A_{11},B_{11},C_{11}+A_{12}R_{21})$
    \hspace*{1.3in} ... solve (\ref{eqn_Syl2}) \\
(c) \> \> \> $R_{22} = SylR(A_{22},B_{22},C_{22}-R_{21}B_{12})$
    \hspace*{1.3in} ... solve (\ref{eqn_Syl3}) \\
(d) \> \> \> $R_{12} = SylR(A_{11},B_{22},C_{12}-R_{11}B_{12}+A_{12}R_{22})$
    \hspace*{.7in} ... solve (\ref{eqn_Syl4}) \\
\> \> end
\end{tabbing}

If the matrix multiplications in $SylR$ were done using the conventional
$O(n^3)$ algorithm, then $SylR$ would perform the same
arithmetic operations (and so make the same rounding errors) 
as a conventional Sylvester solver, just in a different order.
For the complexity analysis of $SylR$ we assume $n$ is a power of 2:
\begin{eqnarray*}
cost(n) & = & cost(n/2) 
                \hspace{1.9in} ... {\rm cost\ of\ line\ (a)} \\
          &   &  + MM(n/2) + cost(n/2) 
                \hspace{1in} ... {\rm cost\ of\ line\ (b)} \\
          &   &  + MM(n/2) + cost(n/2) 
                \hspace{1in} ... {\rm cost\ of\ line\ (c)} \\
          &   &  + 2 \cdot MM(n/2) + cost(n/2) 
                \hspace{.85in} ... {\rm cost\ of\ line\ (d)} \\
        & = & 4 \cdot cost(n/2) + O(n^{\omega + \eta}) \\
        & = & O(n^{\omega + \eta})  
                \hspace{2in} ... {\rm as\ long\ as\ } \omega + \eta > 2
\end{eqnarray*}

Proving logarithmic stability will depend on each subproblem having
a condition number bounded by the condition number of the original
problem. Since the Sylvester equation is really the triangular
linear system $(I \otimes A - B^T \otimes I) \cdot \VEC (R) = - \VEC (C)$, where
$\otimes$ is the Kronecker product and $\VEC (R)$ is a vector of
the columns of $R$ stacked atop one another from left to right,
the condition number of the Sylvester equation is taken to be the
condition number of this linear system. This in turn is governed
by the smallest singular value of the matrix, which is denoted $\sep (A,B)$:
\[
\sep (A,B) \equiv \sigma_{\min} (I \otimes A - B^T \otimes I) 
   = \min_{\|R\|_F = 1} \| AR-RB \|_F
\]
where $\| \cdot \|_F$ is the Frobenius norm \cite{varah}.
Just as in Section~\ref{sec_RTriInv}, where each subproblem
involved inverting a diagonal block of the triangular
matrix, and so had a condition number bounded by the
original problem, here each subproblem also satisfies
\begin{equation}
\label{eqn_SepBound}
\sep (A_{ii},B_{jj}) \geq \sep (A,B) \; \; .
\end{equation}

This lets us prove that this algorithm is
logarithmically stable as follows.
Similar to before, we use the induction hypothesis that $err(n')$ bounds the
error in the solution of any smaller Sylvester equation of dimension
$n'$ encountered during the algorithm (including errors in computing the
right-hand-side). We use the fact that changing the right-hand-side of a
Sylvester equation by a matrix bounded in norm by $x$ can change the solution
in norm by $x/\sep (A,B)$, as well as error bound~(\ref{eqn_MulBnd}):
\begin{align*}
err(R_{21},n/2) & \leq err(n/2)  \\
err(R_{11},n/2) & \leq  
    err(n/2) + \frac{1}{\sep (A,B)} ( 
    \varepsilon \| C_{11} \| 
  + \| A_{12} \| \cdot err (R_{21},n/2)
  + \mu(n/2) \varepsilon \| A_{12} \| \cdot \|R_{21} \| ) \\
                & \leq 
    err(n/2) + \frac{1}{\sep (A,B)} ( 
    \varepsilon \| C \| 
  + \| A \| \cdot err (n/2)
  + \mu(n/2) \varepsilon \| A \| \cdot \|R\| ) \\
err(R_{22},n/2) & \leq  
    err(n/2) + \frac{1}{\sep (A,B)} ( 
    \varepsilon \| C_{22} \| 
  + \| B_{12} \| \cdot err (R_{21},n/2)
  + \mu(n/2) \varepsilon \| B_{12} \| \cdot \|R_{21} \| ) \\
                & \leq 
    err(n/2) + \frac{1}{\sep (A,B)} ( 
    \varepsilon \| C \| 
  + \| B \| \cdot err (n/2)
  + \mu(n/2) \varepsilon \| B \| \cdot \|R\| ) \\
err(R_{12},n/2) & \leq  
    err(n/2) + \frac{1}{\sep (A,B)} ( 
    \varepsilon \| C_{12} \| 
  + \| B_{12} \| \cdot err (R_{11},n/2)
  + \mu(n/2) \varepsilon \| B_{12} \| \cdot \|R_{11} \| \\
  & \hspace*{.5in}
  + \| A_{12} \| \cdot err (R_{22},n/2)
  + \mu(n/2) \varepsilon \| A_{12} \| \cdot \|R_{22} \| ) \\
                & \leq 
    err(n/2) + \frac{1}{\sep (A,B)} ( 
    \varepsilon \| C \| 
  + (\| A \| + \| B \| ) \cdot err (n/2)
  + \mu(n/2) \varepsilon ( \| A \| + \| B \| ) \cdot \|R\| ) \\
err(R,n) & \leq  err(R_{11},n/2) + err(R_{12},n/2) + 
                 err(R_{21},n/2) + err(R_{22},n/2) \\
         & \leq  \left( 4 + 2 \frac{\|A\| + \|B\|}{\sep (A,B)} \right) 
                  \cdot err(n/2) + 
                 \frac{\varepsilon}{\sep (A,B)} 
                 ( 3 \|C\| + 2 \mu (n/2) (\|A\| + \|B\|) \|R\| ) \\
\end{align*}
This yields the recurrence
\begin{eqnarray}
\label{eqn_SylError}
err(n) &   =  &   \left( 4 + 2 \frac{\|A\| + \|B\|}{\sep (A,B)} \right) \cdot
               err(n/2) + 
               \frac{\varepsilon}{\sep (A,B)} 
               ( 3 \|C\| + 2 \mu (n/2) (\|A\| + \|B\|) \|R\| ) \nonumber \\
       & \leq &  O( \frac{ n \varepsilon}{\sep (A,B)} 
               ( \|C\| + \mu (n/2) (\|A\| + \|B\|) \|R\| )
               (2 + \frac{\|A\| + \|B\|}{\sep (A,B)})^{\log_2 n} ) \nonumber \\
       & \leq &  O \left( n^{1+\log_2 3} \mu(n/2) \varepsilon \|R\|
                 \left( \frac{\|A\| + \|B\|}{\sep (A,B)} \right)^{1 + \log_2 n}
                 \right) \\
\end{eqnarray}
Since the error bound of a conventional algorithm for the
Sylvester equation is 
$O( \varepsilon \|R\| \frac{\|A\| + \|B\|}{\sep (A,B)} )$,
we see our new algorithm is logarithmically stable.

A similar approach works for the generalized Sylvester equation,
which we just sketch. Equation~(\ref{eqn_Syl}) becomes
\begin{eqnarray*}
\bmat{cc} A_{11} & A_{12} \\ 0      & A_{22} \emat \cdot
\bmat{cc} R_{11} & R_{12} \\ R_{21} & R_{22} \emat -
\bmat{cc} L_{11} & L_{12} \\ L_{21} & L_{22} \emat \cdot
\bmat{cc} B_{11} & B_{12} \\ 0      & B_{22} \emat 
& = &
- \bmat{cc} C_{11} & C_{12} \\ C_{21} & C_{22} \emat \\
\bmat{cc} \bar{A}_{11} & \bar{A}_{12} \\ 0      & \bar{A}_{22} \emat \cdot
\bmat{cc} R_{11} & R_{12} \\ R_{21} & R_{22} \emat -
\bmat{cc} L_{11} & L_{12} \\ L_{21} & L_{22} \emat \cdot
\bmat{cc} \bar{B}_{11} & \bar{B}_{12} \\ 0      & \bar{B}_{22} \emat 
& = &
- \bmat{cc} \bar{C}_{11} & \bar{C}_{12} \\ \bar{C}_{21} & \bar{C}_{22} \emat 
\; \; .
\end{eqnarray*}
Multiplying these out leads to four generalized Sylvester equations, one
for each subblock, which are solved consecutively and recursively as before.

\subsection{Computing (generalized) Eigenvectors}
\label{sec_Eigenvectors}

Given a matrix in Schur form, i.e. triangular, each eigenvector
may be computed in $O(n^2)$ operations by solving a triangular
system of equations. But this would cost $O(n^3)$ to compute 
all $n$ eigenvectors, so we need a different approach.
Rewriting equation~(\ref{eqn_InvariantSS}) as
\begin{equation}
\label{eqn_BlockDiag}
T \equiv
\bmat{cc} A & C \\ 0 & B \emat
=
\bmat{cc} I & R \\ 0 & I \emat
\cdot
\bmat{cc} A & 0 \\ 0 & B \emat
\cdot
\bmat{cc} I & R \\ 0 & I \emat^{-1}
\end{equation}
we see that we have reduced the problem of computing
an eigendecomposition of $T$ to finding
eigendecompositions of $A = V_A \Lambda_A V_A^{-1}$ and 
$B = V_B \Lambda_B V_B^{-1}$ separately, and then 
multiplying \linebreak
$\bmat{cc} I & R \\ 0 & I \emat \cdot \bmat{cc} V_A & 0 \\ 0 & V_B \emat$
to get the eigenvector matrix of $T$. This leads to the following
divide-and-conquer algorithm, where as before we assume
all submatrices are square, of the same power-of-2 dimensions:

\begin{tabbing}
jnk \= jnk \= jnk \= jnk \= jnk \kill
\> function $V = EVecR(T)$ ... all matrices are $n$-by-$n$ \\
\> \> if ($n=1$) then \\
\> \> \> $V = 1$ \\
\> \> else ... use notation from equation~(\ref{eqn_BlockDiag}) \\
(a) \> \> \> $R = SylR(A,B,C)$ \\
(b) \> \> \> $V_A = EVecR(A)$ \\
(c) \> \> \> $V_B = EVecR(B)$ \\
(d) \> \> \> $V = \bmat{cc} V_A & R \cdot V_B \\ 0 & V_B \emat$ \\
(e) \> \> \> for $i=n/2+1$ to $n$, $V(:,i) = V(:,i) / \| V(:,i) \|_2$, end for \\
\> \> end if
\end{tabbing}

For the complexity analysis we assume as before that $n$ is a power of 2:
yielding
\begin{eqnarray*}
cost(n) & = & O(n^{\omega + \eta})
                \hspace{1.4in} ... {\rm cost\ of\ line\ (a)} \\
          &   &  + 2 \cdot cost(n/2) 
                \hspace{1.1in} ... {\rm cost\ of\ lines\ (b)\ and\ (c)} \\
          &   &  + MM(n/2) 
                \hspace{1.25in} ... {\rm cost\ of\ line\ (d)} \\
          &   &  + O(n^2)
                \hspace{1.5in} ... {\rm cost\ of\ line\ (e)} \\
        & = & 2 \cdot cost(n/2) + O(n^{\omega + \eta}) \\
        & = & O(n^{\omega + \eta}) 
                \hspace*{1.4in} {\rm as\ desired.}
\end{eqnarray*}

In general, each eigenvector has a different condition number, 
and ideally should be computed with a corresponding accuracy. 
If we compute each eigenvector separately using the conventional 
algorithm in $O(n^2)$ operations, this will be the case, but
it would take $O(n^3)$ operations to compute all the eigenvectors
this way. 

Unfortunately, $EVecR$ cannot guarantee this accuracy, for the 
following reason. If the first splitting of the spectrum is
ill-conditioned, i.e. $\sep (A,B)$ is tiny, then this will
affect the accuracy of all subsequently computed
right eigenvectors for $B$, through the multiplication $R \cdot V_B$.
Indeed, multiplication by $R$ can cause an error in any 
eigenvector of $B$ to affect any other eigenvector.
(It would also affect the left (row) eigenvectors
$[ V_A^{-1}, -V_A^{-1} \cdot R ]$ for $A$, if we computed them.)
This is true even if some of $B$'s eigenvectors are much better
conditioned. Similarly, further splittings within $A$ (resp. $B$) 
will effect the accuracy of other right eigenvectors of $A$
(resp. $B$), but not of $B$ (resp. $A$),


To simplify matters, we will derive one error bound that works
for all eigenvectors, which may overestimate the error for some.
$\| \cdot \|$ will denote the 2-norm.
We will do this by using a lower bound $\underline{s}>0$ for all
the values of $\sep (A,B)$ encountered during the algorithm.
Thus $\underline{s}$ could be taken to be the minimum value 
of all the $\sep (A,B)$ values themselves, or bounded below
by $\min_{1 \leq i < n} \sep (T(1:i,1:i),T(i+1:n,i+1:n))$,
which follows from inequality~(\ref{eqn_SepBound}).
We can also use $\| R \| \leq \| T \| / \underline{s}$ as
a bound on the norm of any $R$ at any stage in the algorithm.
We use this to simplify bound~(\ref{eqn_SylError}) on the
error of solving the Sylvester equation, yielding
\begin{equation}
\label{eqn_SylError2}
O \left( n^{c} \varepsilon 
\left( \frac{\|T\|}{\underline{s}} \right)^{2 + \log_2 n} \right) 
\end{equation}
for a modest constant $c$.
Using the induction hypothesis that $err(n')$ bounds the
error in the computed $n'$-by-$n'$ eigenvector matrix at
any stage in the algorithm, as well as bound~(\ref{eqn_MulBnd}),
we get
\begin{align*}
err(R, n/2) & =
O \left( n^{c} \varepsilon 
\left( \frac{\|T\|}{\underline{s}} \right)^{2 + \log_2 n} \right) \\
err(V_A, n/2) & \leq err(n/2) \\
err(V_B, n/2) & \leq err(n/2) \\
err(V,n)      & \leq err(V_A, n/2) + err(V_B, n/2) + 
                     err(R, n/2) \|V_B\| + \|R\| err(V_B) + 
                     \mu(n/2) \varepsilon \|R\| \cdot \|V_A\| \\
              & \leq (2 + \|R\|) err(n/2) + \sqrt{n} \cdot err(R,n/2)
              + \sqrt{n} \mu (n/2) \varepsilon \|R\| \\
              & \leq (2 + \|R\|) err(n/2) + 
                     O \left( n^{c'} \varepsilon 
                     \left( \frac{\|T\|}{\underline{s}} \right)^{2 + \log_2 n} \right) \\
\end{align*}
for another modest constant $c'$. Step (e) of $EVecR$, which makes each column have
unit norm, will decrease the absolute error in large columns, but 
we omit this effect from our bounds, which are normwise in nature, and so get
a larger upper bound.
Bounding $2 + \|R\| \leq 3 \frac{\|T\|}{\underline{s}}$,
and changing variables from $n=2^m$ to $m$ and $err(n)$ to $\bar{err}(m)$, we get
\[
\bar{err}(m) \leq 3 \frac{\|T\|}{\underline{s}} \cdot \bar{err}(m-1) +
                     O \left( \varepsilon 
                     \left( \frac{2^{c'}\|T\|}{\underline{s}} \right)^{2 + m} \right) 
\]
Finally, setting $f(m) = \bar{err}(m)/ (3 \frac{\|T\|}{\underline{s}})^m$, we get
a simple geometric sum
\[
f(m) \leq f(m-1) + O \left( \left( \frac{2^{c'}}{3} \right)^m \varepsilon 
                     \left( \frac{\|T\|}{\underline{s}} \right)^{2} \right) 
\]
which, since $2^{c'}>3$, leads to our final bound
\[
err(n) = O \left( 
n^{c'} \varepsilon \left( \frac{\|T\|}{\underline{s}} \right)^{2 + \log_2 n}
\right) \; \; ,
\]
demonstrating a form of logarithmic stability.



\section{Conclusions}

We have shown that nearly all standard dense linear algebra operations
(LU decomposition, QR decomposition, matrix inversion, 
linear equation solving,
solving least squares problems, computing the (generalized) Schur form,
computing the SVD, and solving (generalized) Sylvester equations)
can be done stably and 
asymptotically as fast 
as the fastest matrix multiplication algorithm that may ever 
exist (whether the matrix multiplication
algorithm is stable or not).  
For all but matrix inversion and solving (generalized) Sylvester equations,
stability means backward stability in
a normwise sense, and we measure complexity by
counting arithmetic operations.

For matrix inversion and solving the Sylvester equation,
stability means forward stability, 
i.e. that the error is bounded in norm by 
$O( \varepsilon \cdot \kappa )$, 
machine epsilon times the appropriate condition number,
just as for a conventional algorithm. 
The conventional matrix inversion algorithm is not
backward stable for the matrix as a whole either,
but requires a different backward error for each column.
The conventional solution of the Sylvester equation
is not backward stable either, and only has a forward error bound.

Also, for matrix inversion and solving the Sylvester equation, we
measure complexity by counting bit operations, to account
for the use of a modest amount of extra precision.
Indeed, we can say that matrix multiplication (stable or not)
in $O(n^{\omega + \eta})$ operations for any $\eta > 0$ 
is possible if and only if forward stable inversion of 
an $n$-by-$n$ matrix $A$ in $O(n^{\omega + \eta})$ operations 
for any $\eta > 0$ is possible.
(See Theorem~\ref{thm_1} for a more careful statement of
how operations are counted.)

All eigenvectors of (generalized) nonsymmetric eigenvalue problems
can also be computed in $O(n^{\omega + \eta})$ bit operations
from the Schur form,
but with a weaker notion of forward stability, where the error
bound for all the eigenvectors, in the worst case, 
depends on the largest condition number of any eigenvector.

Finally, we note several possible practical implications
of our algorithms.
Several of the recursive algorithms
we used (QRR and LUR) were originally invented for their
superior memory locality properties 
\cite{ElmrothGustavson2000,Toledo97}, and
the same property is likely to hold for our
new recursive algorithms as well. 
The divide-and-conquer nature of these algorithms
also naturally creates parallelism that could be
exploited on various architectures. 
Even if one is not using an asymptotically faster
matrix-multiplication algorithm, these algorithms
could be advantageous on platforms that perform
matrix multiplication much faster than other 
basic linear algebra subroutines.

\section{Acknowledgments}

The authors would like to thank Robert Kleinberg for
many useful discussions, as well as an anonymous referee for
many detailed and useful comments.


\bibliographystyle{plain}
\bibliography{linalg,biblio1,matmulJD,matrix2}

\end{document}